\documentclass[final]{siamltex}
\usepackage{xcolor}
\usepackage{amsmath}
\usepackage{amssymb}
\usepackage{mathrsfs}
\usepackage{graphicx}
\usepackage{hyperref}
\usepackage{bm}
\usepackage{ucs}
\usepackage[utf8x]{inputenc}
\usepackage{wrapfig}
\usepackage{color}
\usepackage{float}
\usepackage{epstopdf}

\usepackage{amsmath,amstext,amssymb,graphicx}
\usepackage{verbatim}
\usepackage{color}
\usepackage{algorithm,algorithmic}

\usepackage{mathtools}

\usepackage{epstopdf}
\usepackage{color}
\usepackage{graphicx}
\usepackage{float}
\usepackage[section]{placeins}

\usepackage{tikz}
\usetikzlibrary{arrows}

\def\Gc{{\cal G}}

\def\norm#1{\|#1\|}

\def\vect#1{\mbox{\boldmath{$#1$}}}

\def\mC{\mathbb{C}}

\def\wG{\widehat{G}}
\def\wU{\widehat{U}}
\def\wV{\widehat{V}}

\def\we{\widehat{e}}

\def\wg{\widehat{g}}

\def\wP{\widehat P}

\newcommand{\Pm}{\vect P}
\newcommand{\Mm}{\vect M}
\newcommand{\Pp}{P_{\!\mbox{\scriptsize{parax}}}}

\newcommand{\eps}{\varepsilon}
\renewcommand{\epsilon}{\eps}
\renewcommand{\leq}{\leqslant}

\def\bfrho{\mbox{\boldmath$\rho$}}

\DeclarePairedDelimiter{\abs}{\lvert}{\rvert}

\begin{document}

\title{Coherent imaging without phases}
\author{Miguel Moscoso, Alexei Novikov and George Papanicolaou}
\maketitle
\begin{abstract}
In this paper we consider narrow band, active array imaging of weak localized scatterers
when only the intensities are recorded at an array with $N$ transducers. We assume that
the medium is homogeneous and, hence,  wave propagation is 
fully coherent. 
This work is an extension of our previous paper~\cite{Novikov14} where 
we showed that using linear combinations of intensity-only measurements, obtained from $N^2$ illuminations,
imaging of localized scatterers can be carried out  efficiently using imaging methods based on the singular value decomposition of 
the time-reversal matrix. 
Here we show the same strategy can be accomplished with only $3N-2$ illuminations, therefore reducing enormously the data acquisition
process. Furthermore, we show that in the paraxial regime one can form the images by using 
six illuminations only.  In particular, this
 paraxial  regime includes 
 Fresnel and Fraunhofer diffraction.  The key point of this work is that if one controls the illuminations, imaging with 
intensity-only can be easily reduced to a  imaging 
with phases and, therefore, one can apply standard imaging techniques. Detailed numerical simulations illustrate the performance of 
the proposed imaging strategy with and without data noise.
\end{abstract}

\section{Introduction}
In many situations it is difficult, or impossible, to measure the phases received at the detectors, only the intensities are available for imaging.
This is the case, for example, in imaging from X-ray sources \cite{Millane90, Harrison93, Pfeiffer06},
 or from optical sources \cite{Walther63,Dainty87,Trebino93}, where one wants to form images from the spectral intensities. This is known as the phase retrieval problem, in which one seeks to reconstruct a complex signal from quadratic 
measurements involving the signal and, possibily, additional {\em a priori} information about it. The literature on the subject ranges from more or less sophisticated experimental setups that use interferometry to retrieve the phases
\cite{Yamaguchi97,Cuche99},  to the use of algorithms that obtain the phases of the signals received at the array and then  form the images \cite{GS72,Fienup82,CMP11,Candes13}.

By far, the most popular methods for reconstructing a signal  
from only the magnitudes of its Fourier coefficients are alternating projection algorithms proposed by Gerchberg and Saxton \cite{GS72} and later improved by Fienup \cite{Fienup82}.
These algorithms use a sequence of efficient projections in the Fourier and the spatial domains to reconstruct the missing phases of the Fourier coefficients that
are consistent with their magnitudes and with the known spatial constrains. However, it is well known that these algorithms do not guarantee convergence, as they may suffer from stagnation of the iterates away from the true solution.  
In particular, they do not work well when prior knowledge about the sought signal is not available or is poor, and they often need a careful usage of the domain constrains (such as support and nonnegativity).

Another class of algorithms were proposed by Chai et. al. in~\cite{CMP11} that framed the problem of array imaging using intensity-only 
measurements as a low-rank matrix completion problem. 
Motivated by the theory in compressed sensing, and in particular by the recent developments in matrix rank minimization theory \cite{Candes09}, 
the authors in~\cite{CMP11}  replace  the  original non linear vector problem 
by a linear matrix problem with a rank-one solution. Since the resulting rank minimization problem is NP-hard, they relaxed it into a convenient convex program that seeks a low rank matrix 
that matches the intensity data by using nuclear norm  minimization.  
A similar approach was used by Candes et. al. in~\cite{Candes13} for diffraction imaging.

The main advantage of the approach in \cite{CMP11,Candes13} is that it guarantees exact recovery under some conditions on the imaging operator and the sparsity of the image. No additional constraints are needed on the image. However,
the fact that a vector problem is replaced by a  matrix one makes the  search of the solution very expensive if the  images are large. Indeed, if the image is discretized using $K$ pixels, then the reconstruction of the unknown image is done in a  space of $K^2$ dimensions. 
This grand optimization problem is not feasible
if, for example, a high resolution of the image is sought and, hence, the number of pixels is large.

The objective of this paper is to propose a new perspective for active array imaging of weak localized scatterers when
only the intensities are available for imaging. We follow the same approach as in \cite{Novikov14}, where we came up with 
a new strategy that guarantees exact recovery and that is efficient for large scale problems.
In \cite{Novikov14}, we showed that imaging with intensity-only measurements can be carried out 
by using an appropriate protocol of illuminations and the polarization identity. This allows us to obtain
the time reversal operator of the imaging system. Once the time reversal operator has been obtained, the images can be formed using its singular value decomposition (SVD). This approach guarantees exact recovery of the image while keeping the computational cost low. However, the proposed strategy required $N^2$ illuminations and,
therefore, the data acquisition process was expensive. 
In this paper we show that, in fact, only $3N-2$ illuminations are needed, making the proposed approach in \cite{Novikov14} feasible in practice. 

We also consider the paraxial regime, which is the appropriate scaling for Fresnel or Fraunhofer diffraction. In this regime, array imaging has the form of a Fourier transform and, hence, the process of recovering an image is the classic phase retrieval problem. In this case, we show that only six illuminations 
are needed to create an image of a flat object. In our case, a flat object consists of a set of point-like scatterers  located at the same range. 
The key point of the approach proposed here is that multiple versions of the intensity 
of the Fourier transform of an 
object, obtained through an appropriate set of illuminations, make
the solution of the inverse problem of phase retrieval unique. We note, however, that the use of 
multiple measurements to resolve phase uniqueness is not new. Redundancy of the data to 
enforce phase uniqueness can be obtained, for example, by using 
random illuminations as in~\cite{Fannjiang12}, or by inserting masks as in~\cite{Candes13}. 
Another very interesting possibility is used in ptychography~\cite{Rodenburg08}, a form of 
coherent diffractive imaging, where the object is stepped through a localized coherent wavefront 
generating a series of diffraction patterns.
The illuminated area at each position overlaps with its neighbors and, thus, redundancy can be 
exploited during iterative phase-retrieval ptychography.


We also carry out numerical simulations that address the limitations of the proposed approach when the phaseless data is noisy.
We find that sensitivity to noise is higher in the Fraunhofer regime when fewer illuminations are needed.
That is the method of six illuminations, which is appropriate when the distance between the scatterers and the 
array is much larger than the wavelength of the signals and much larger than the linear dimensions of the array and the IW, is very sensitive to 
noise. This indicates that there is a trade-off between using a limited number of illuminations with phaseless data and the level of noise in the data.

The paper is organized as follows. In Section \ref{sec:model}, we formulate the active array imaging problem using intensity-only
measurements. In Section \ref{sec:illum} we describe the illumination strategy and the imaging approach proposed in the paper.
Section \ref{sec:numerics} contains our numerical simulations. The conclusions of this work are in Section \ref{sec:conclusions}.

\section{Active array imaging}
\label{sec:model}
In active array imaging we seek  to determine the location and reflectivities of a few reflectors by sending probing signals from an array and
recording the backscattered signals. The active array ${\cal A}$ consists of $N$ transducers placed at distance $h$ between them. 
These transducers emit spherical wave signals from positions $\vec{\vect x}_s\in{\cal A}$ and record the echoes with receivers at positions $\vec{\vect x}_r\in{\cal A}$.  
In this paper, we consider narrow-band array imaging of sparse images consisting of a few point-like scatterers in a homogeneous medium. By point-like scatterers we mean very small scatterers compared to the wavelength of the probing signals. 
We assume that  multiple scattering between the scatterers is negligible. For imaging problems with
multiple scattering see \cite{Chai14}.
\begin{figure}[htbp]
\begin{center}
\includegraphics[scale=0.9]{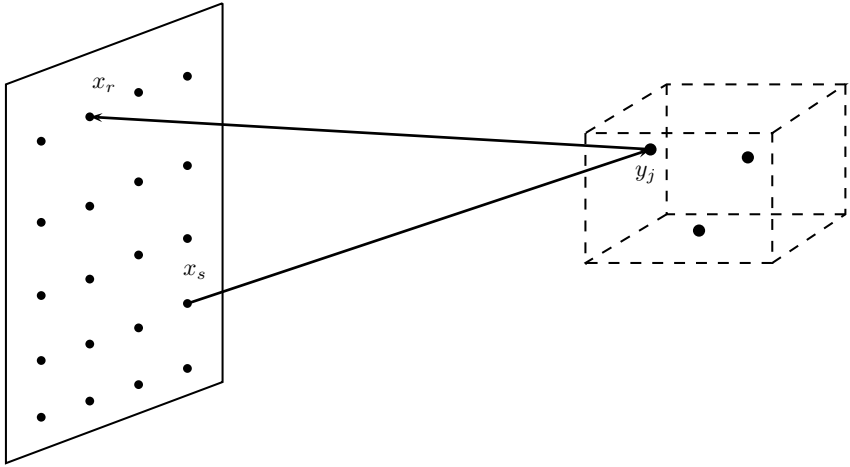}
\caption{A general setup of array imaging problem}
\label{cs_3d_illustration}
\end{center}
\end{figure}

A typical configuration of the imaging setup of the array imaging problem is given in Figure~\ref{cs_3d_illustration}. 
The active array has $N$ transducers at positions $\vec{\vect x}_s$
located on the plane $z=0$, so $\vec{\vect x}_s=(\vect x_s,0)$, $s=1,\cdots,N$. 
There are $M$ point-like scatterers in the image window (IW)
that is discretized using a uniform grid of $K$ points 
$\vec{\vect y}_j$, $j=1,\ldots,K$. The images are sparse because $K\gg M$.
The scatterers with reflectivities $\alpha_j\in\mC$ are
located at  positions $\vec{\vect \xi}_{1},\ldots,\vec{\vect \xi}_{M}$, which we assume coincide with one of these $K$ grid points so
$\{\vec{\vect \xi}_{1},\ldots,\vec{\vect \xi}_{M} \} \subset \{ \vec{\vect y}_1, \dots \vec{\vect y}_K \}$. 
For a study on imaging with $\ell_1$ optimization when the scatterers lie off-grid we refer to \cite{Fannjiang12b,Borcea15}. 
If the scatterers are far apart or the reflectivities are small, the interaction between them is weak and  the Born approximation is applicable. 
In this case, the response at $\vec{\vect x}_r$ due to a narrow-band pulse of angular frequency $\omega$ sent from $\vec{\vect x}_s$ and reflected by the $M$ scatterers is given by
\begin{equation}
\label{response}
P(\vec{\vect x}_r,\vec{\vect x}_s)=\sum_{j=1}^M\alpha_j G (\vec{\vect x}_r,\vec{\vect \xi}_j) G (\vec{\vect \xi}_{j},\vec{\vect
x}_s)\, ,
\end{equation}
where
\begin{equation}
\label{homo_green}
G(\vec{\vect x},\vec{\vect y}) = \frac{\exp\{i \kappa \abs{\vec{\vect x} - \vec{\vect y}}\}}{4\pi \abs{\vec{\vect x} - \vec{\vect y}}} 
\end{equation}
is the Green's function that characterizes 
wave propagation from $\vec{\vect x}$ to $\vec{\vect y}$ in a 3-dimensional homogeneous medium. In \eqref{homo_green}, $\kappa=\omega/c_0$ 
is the wavenumber and $c_0$ is the wave speed in the medium. 
The nonlinear problem that includes multiple scattering between the scatterers is considered in \cite{Chai14}.

To write the data received on the array in a more compact form, we define the {\em Green's function vector} 
$\vect g (\vec{\vect y})$ at location $\vec{\vect y}$ in IW as 
\begin{equation}
\label{GreenFuncVec}
\vect g (\vec{\vect y})=[ G(\vec{\vect x}_{1},\vec{\vect y}),\cdots,G(\vec{\vect x}_{N},\vec{\vect y})]^t\, ,
\end{equation}
where $.^t$ denotes  transpose. This vector represents the signal received at the array due to a point source at $\vec{\vect y}$. 
It can also be interpreted as the illumination vector of the array targeting the position $\vec{\vect y}$. 
We also introduce the true {\it reflectivity vector}
$\vect\rho=[\rho_{1},\ldots,\rho_{K}]^t\in\mC^K$ such that
\begin{equation}
\label{eq:rho}
\rho_{k}=\sum_{j=1}^M\alpha_j\delta_{\vect \xi_{j},\vect y_k},\,\, k=1,\ldots,K,
\end{equation}
where $\delta_{\cdot,\cdot}$ is the classical Kronecker delta. With this notation
we can write the response matrix  as a sum of rank-one matrices so
\begin{equation}
\label{responsematrix}
\Pm  \equiv [ P(\vec{\vect x}_r,\vec{\vect x}_s)]_{r,s=1}^{N}=\sum_{j=1}^M\alpha_j\vect g(\vect \xi_{j})\vect g^t(\vect \xi_{j})
=\sum_{j=1}^K\rho_{j}\vect g (\vect \xi_{j})\vect g^t(\vect \xi_{j}).
\end{equation}
The entry $P(\vec{\vect x}_r,\vec{\vect x}_s)$ of this matrix denotes the signal received at $\vec{\vect x}_r$ due to a
signal sent from $\vec{\vect x}_s$.
Using \eqref{GreenFuncVec}, we also define the $N\times K$ sensing matrix 
\begin{equation}\label{sensingmatrix}
\vect{\Gc}=[\vect g(\vect y_1)\,\cdots\,\vect g (\vect y_K)] \,,
\end{equation}
whose column vectors are the signals received at the array due to point sources at the grid points. Thus,
$\vect{\Gc}$ maps a distribution of sources in the IW to the data received on the array. 
Using \eqref{sensingmatrix}, we write \eqref{responsematrix} in matrix form as
\begin{equation}\label{responsematrix1}
\Pm =\vect\Gc \diag(\bfrho )\vect\Gc^t.
\end{equation}
The full response matrix  $\Pm$ represents a linear transformation from the illumination space $\mC^N$ to the data space $\mC^N$. 
Indeed, consider an illumination vector 
\[ \vect f =[f_{1},\ldots, f_{N}]^t\, ,  
\]
whose components are the signals $f_{1},\ldots,f_{N}$ 
sent from each of the $N$ transducers in the array. Then, $\vect\Gc^t\vect f$
gives the signals at each grid point of the IW. These signals are reflected by the scatterers on the 
grid that have reflectivities given by the vector $\bfrho$, and then they are backpropagated to the array by the matrix $\vect\Gc$.
All the  information available for imaging, including phases, is contained in the full response matrix  $\Pm$. If
sources and receivers are located at the same positions, then it is symmetric (due to Lorentz reciprocity) but not hermitian.

Given a set of illuminations $\{\vect f^{(j)} \}_{j=1,2,\dots}$, the usual imaging problem is to determine the location and reflectivities of the scatterers 
from  the data 
\begin{equation}\label{data}
\vect b^{(j)}=\Pm \vect f^{(j)}\,\, , \quad j=1,2,\dots \, ,
\end{equation}
recorded at the array, including phases. 
If phase information is not available because only the intensities $\beta^{(j)}_i = \abs{ b^{(j)}_i }^2$ of the signals, $i=1,\dots,N$, are recorded at the array, the imaging problem is to determine the location and reflectivities of the scatterers from  the absolute values of each component in \eqref{data}, i.e., from the intensity vectors
\begin{equation}\label{dataI}
\vect \beta^{(j)}=\diag \left( \left(\vect P \vect f^{(j)} \right) \left( \vect \Pm \vect f^{(j)} \right)^* \right)\,\, , \quad j=1,2,\dots \, ,
\end{equation}
where the superscript $*$ denotes conjugate transpose. 

\section{Coherent illuminations}
\label{sec:illum}
The essential point in active array imaging of scatterers is that it is possible to control the illuminations to form the images. 
The main objective of this paper is to show that if one controls the illuminations the imaging problem with intensity-only
can be reduced to one in which both intensities and phases are available and, therefore,
algorithms that use the full dataset can be used to form the images. In this section, we show how to obtain the missing 
phases of the signals recorded at the array from linear combinations of  the intensities of these signals. We first consider the general case in which 
sources and receivers are not placed at the same locations and, therefore, the response matrix of the imaging system is not symmetric. Then, we consider the 
typical imaging setup in which sources and receivers are located at the same positions and, therefore,
the response matrix is symmetric. We finally consider the important case of Fresnel and Fraunhofer diffraction, where wave propagation can be modeled by the paraxial approximation.
\subsection{General case}
\label{sec:illum_general}
In~\cite{Novikov14}, we proposed a novel imaging strategy for the case in which only data of the form \eqref{dataI}
are recorded at the array. The main idea behind that approach is to use
 the {\em time reversal operator} $\Mm= \Pm^* \Pm$ for imaging, which can be obtained from
intensity measurements using an appropriate illumination strategy and the polarization identity. Here
 we show how to obtain $\Mm$ for imaging using a much more efficient illumination strategy. We do not assume that the response matrix $\Pm$ is symmetric and, therefore, sources and receivers do not need to be located at the same positions.

In~\cite{Novikov14} we showed how to obtain all the entries of the matrix $\Mm$ using $O(N^2)$ illuminations, where
$N$ is the number of sources in the array. That method 
 uses only the total power received by the entire array, and 
 the following polarization identities
 \begin{equation}\label{polarization_re}
\mbox{Re}(M_{ij}(\omega))= 
\frac{1}{2} \left( \norm{\vect\wP(\omega) \vect\we_{i+j}}^2  - 
\norm{\vect\wP(\omega) \vect\we_i}^2 - \norm{\vect\wP(\omega) \vect\we_j}^2\right)
\end{equation}
\begin{equation}\label{polarization_im}
\mbox{Im}(M_{ij}(\omega))= 
 \frac{1}{2} \left( \norm{\vect\wP(\omega) \vect\we_{i-{\bf i} j}}^2  -
  \norm{\vect\wP(\omega) \vect\we_i}^2 - \norm{\vect\wP(\omega) \vect\we_j}^2\right), 
  \end{equation}
where $\vect\we_{i} = [0, 0, . . . , 1, 0, . . . , 0]^T$ is the illumination vector whose entries are all zero
except the i-th entry which is $1$. In \eqref{polarization_re}-\eqref{polarization_im},
\[
\vect\we_{i+j}= \vect\we_{i}+\vect\we_{j},  \quad \vect\we_{i-{\bf i} j}= \vect\we_{i}-{\bf i} \vect\we_{j}.
\]

If instead of using  the total power measured at the array as in~\cite{Novikov14}, we use  the intensities of the signals measured at every receiver 
separately, then we can obtain $\Mm$ with significantly fewer $O(N)$ illuminations. 
Observe that each entry of the {\em time reversal operator} $\Mm= \Pm^* \Pm$ can be written as
\[
m_{ij} =\sum_{k=1}^N p_{k i} \bar{p}_{k j}.
\]
In order to recover $p_{k i} \bar{p}_{k j}$ it suffices to know the amplitudes of the signals $|p_{k i}|$, $|p_{k j}|$, $k=1,\dots, N$, and to 
find the {\it phase differences} $ \arg p_{k i} - \arg p_{k j}$, $k=1,\dots, N$. Here, $ \arg p_{k i}$ denotes the unrecorded phase of signal measured at the  
i-th receiver when the  k-th source emmits a signal. The amplitudes are recorded using  $N$ illuminations $\vect\we_{i}$, $i=1,2,\dots,N$, 
and the phase differences can be recovered as follows.
Since  
\[
  \arg p_{k i} - \arg p_{k j}=  (\arg p_{k 1} - \arg p_{k j}) - ( \arg p_{k 1} - \arg p_{k i}),
\]
it suffices to find the  phase differences $\arg p_{k 1} - \arg p_{k j}$ for $j=2,\dots, N$, which means that only the phase differences between 
the first column of  $\Pm$ and all the other colums are needed. If all $p_{k 1} \neq 0$, these
phase differences can be found from the $2 N-2$ illuminations $\vect\we_{1+j}, \vect\we_{1-{\bf i} j}$, $j=2,\dots, N$,
and the polarization identities~\eqref{polarization_re},~\eqref{polarization_im}.  When the image is sparse, the assumption $p_{k 1} \neq 0$
 is not restrictive because of the uncertainty principle~\cite{DS}. 
 

\subsection{Symmetric case}
\label{sec:illum_symmetric}
When sources and receivers are located at the same positions, $\Pm$ is symmetric. In this case, we can obtain the full response 
matrix $\Pm$, up to a global phase, from intensity measurements using also $O(N)$ illuminations. 
Indeed, once we find   $\arg p_{k i} - \arg p_{k j}$ for all $i,j,k$, we also know 
$\arg p_{i k} - \arg p_{j k}$ for all $i,j,k$.  Note that in the general case we did not recover $\Pm$ itself, as we did not need it to image the scatterers. 
Indeed, in the general case where $\Pm$ is not symmetric we use $\Mm$ instead, as these two matrices have the same right singular vectors
and, therefore, the images can be formed using SVD-based methods.

We point out that by using $O(s)$ illuminations we can obtain $O(s^2)$ distinct entries of the matrix $\Pm$ up to a global phase. 
The strategy is as follows. Suppose we want 
to recover the $s \times s$ upper left corner of $\Pm$, with $s=3$. 
First, we illuminate with the illumination vector $\vect\we_{1}$, and we set $p_{11}=|p_{11}|$. This fixes the global phase.  Then, 
we use the three illumination vectors  $\vect\we_{2}$,  
$\vect\we_{1}+\vect\we_{2}$, and $\vect\we_{1}-{\bf i} \vect\we_{2}$. Apart from recording $|p_{k 2}|$, $k=1,2,\dots,N$, we also recover 
the relative phases $\arg p_{k 1} - \arg p_{k 2}$, $k=1,2,\dots,N$, using the polarization identity. In particular, since the phase  $\arg p_{11}=0$, we recover 
$p_{1 2}$ with its phase and, since $\Pm$ is symmetric, $p_{2 1}$ is also recovered with its phase. Furthermore, since the relative phase 
$\arg p_{2 1} - \arg p_{2 2}$ is then recovered, we obtain $p_{2 2}$ with its phase as well.  Next, we use the illumination vectors  $\vect\we_{3}$,  
$\vect\we_{1}+\vect\we_{3}$, and $\vect\we_{1}-{\bf i} \vect\we_{3}$. Since  $\arg p_{1 1} - \arg p_{1 3} = - \arg p_{1 3}$ , we recover 
$ p_{1 3}$ with its phase. Again, using the symmetry of $\Pm$ we also obtain $p_{3 1}$. Since we know at this point $|p_{3 2}|$, $| p_{3  3}|$,
$\arg p_{3 1} - \arg p_{3 2}$, and $\arg p_{3 1} - \arg p_{3 3}$, the entire $3 \times 3$ upper left corner of $\Pm$ is recovered. Proceeding
 analogously, we can obtain the $s \times s$ upper left corner of $\Pm$ after $3s-2$ illuminations (up to a global phase).

\subsection{Paraxial Approximation}
\label{sec:illum_six}

Here, we further specialize the illumination strategy for the setup in which the distance between the scatterers and the array is much larger 
than the wavelength of the signals used to probe the medium and much larger than the linear dimensions of the array and the IW. In this case, 
wave propagation can be modeled by the paraxial approximation which is relevant for (near-field) Fresnel and (far-field) Fraunhofer diffraction. 
We assume that  all the scatterers are at the  same range from the array so the images are flat, and that sources and receivers are located at the 
same positions so the full response matrix $\Pm$ is symmetric. In this case, $\Pm$ has, up to pre-factor and post-factor diagonal matrices, 
a simple structure.
More specifically, it can be written as
\begin{equation}
 \vect P =\vect D_{\mbox{\scriptsize receiver}} \vect H \vect D_{\mbox{\scriptsize source}}\, ,
\label{eq:Phankel}
\end{equation}
where $\vect D_{\mbox{\scriptsize receiver}}$ and $\vect D_{\mbox{\scriptsize source}}$ are matrices that only depend on the geometrical setup, and 
$\vect H$ is a matrix with a simple form. In the one-dimensional problem in which we seek 
to image a set of scatterers on a line using a linear array, $\vect H$ is approximately Hankel, i.e., it is a symmetric matrix that is approximately 
constant across the anti-diagonals. 
This special structure of the matrix $\Pm$  allows us to further reduce the number of illuminations to six in the one-dimensional setup, and to twelve in
the two-dimensional setup, independently of the size of the array. 

In the Appendix~\ref{seq:app} we show that when the IW is far enough from the array
the full-response matrix~\eqref{response} 
can be approximated by the paraxial model
\begin{equation}\everymath{\displaystyle}
\Pp(\vec{\vect x}_r,\vec{\vect x}_s)  
= C_r C_s \sum_{j=1}^K \tilde{\rho}_j e^{- i \kappa \frac{\langle \vect x_s + \vect x_r,  
\vect y_j \rangle}{L}}  \,,
\label{eq:fresnel}
\end{equation}
where $\vec{\vect x}_s=({\vect x}_s,0)$ and $\vec{\vect x}_r=({\vect x}_r,0)$ are the positions of the sources and the receivers, respectively, $\vec{\vect y}_j=({\vect y}_j,L), j=1,\dots,K$ are the grid points of the discretized IW, and 
\begin{equation}
\tilde{\rho}_j = {\rho}_j e^{i \kappa \left| \vect y_j \right|^2/2L},  j=1,\dots,K,
\label{eq:m_rho}
\end{equation}
are distorted reflectivities that only modify the phases of the original reflectivitities by a factor that depends on their (unknown) cross-range  positions 
$\vect y_j$, and on the distance $L$ between the array and the IW, which we assume to be known. 
In Eq. \eqref{eq:fresnel},
\begin{equation}\label{eq:factor}
 C_t =C(\vect x_t) = \frac{e^{i  \kappa L}e^{i \kappa \left| \vect x_t \right|^2 /2L}}{4\pi L},~t=r,s,
\end{equation}
are known factors that only depend on the geometrical setup (the layout of the transducers and the distance $L$). 
Then, the approximated response matrix~\eqref{eq:fresnel} can be written as in \eqref{eq:Phankel}
with  diagonal matrices
\[
\vect D_{\mbox{\scriptsize source}}= \vect D_{\mbox{\scriptsize receiver}}  \equiv [D_{r,s} ]_{r,s=1}^{N},  D_{r,s}= \delta_{r,s}C_r.
\]

Equation~\eqref{eq:fresnel} shows that, in the paraxial regime, array imaging with a single illumination
has the form of a Fourier transform. 
Through the introduction of the distorted reflectivities \eqref{eq:m_rho}, it includes both the near-field Fresnel 
regime, for which the Fresnel number
\begin{equation}
F = \frac{a^2}{\lambda L}\, 
\label{eq:fresnelnumber}
\end{equation}
is $O(1)$, and the far-field Fraunhofer regime, for which $F\ll 1$. In the Fraunhofer regime the extra factor in the distorted reflectivities 
(\ref{eq:m_rho}) approaches to $1$, which means that the propagating fronts, when viewed from the IW, are not spherical but planar. 



Because the structure of the full-response matrix is easier to describe in the one-dimensional imaging setup, we consider now a linear array 
on the line $\vec{\vect x} =(x,0,0)$ and a  set of scatterers on the line $\vec{\vect x}' =(x',0,L)$, see Fig. \ref{fig:schematic} (b).  
Assume that the paraxial approximation holds so $P(x_r, x_s) =  \Pp(x_r,x_s)$,   
and consider the processed data 
\begin{equation}
 H(x_r,x_s) = (C_rC_s)^{-1} P(x_r,x_s)=\sum_{j=1}^K \tilde{\rho}_j e^{-i \kappa (x_r +x_s)x'_j/L}\, ,
\label{eq:fullresponseprocessed}
\end{equation}
where $x'_j$ are the grid points of the image, $j=1,\dots,K$. 
Then, for a fixed configuration of the scatterers the entries in~\eqref{eq:fullresponseprocessed} only 
depend on $x_r+x_s$. This shows that if the IW is far enough from the array, the data received on the array only depends on the sum $r+s$ and, 
therefore,  the matrix 
\begin{equation}
\vect H  \equiv [H_{r,s} ]_{r,s=0}^{N-1} 
= [\Xi_{r,s} ]_{r,s=0}^{N-1} 
= \sum_{j=1}^{K} \tilde{\rho}_{j} e^{-i \Lambda \frac{(r+s)j}{N}}\, ,
\quad s,r=0,\dots,N-1\,, 
\label{eq:ptransform}
\end{equation}
has constant skew-diagonals, i.e.,
 $\vect H$ is a Hankel matrix. In \eqref{eq:ptransform}, we have defined
the parameter 
\begin{equation}
 \Lambda = \frac{\kappa}{L}\frac{ab}{K}\, ,
\label{eq:lambda}
\end{equation}
where $a$ and $b$ are the lengths of the array and the IW, respectively. Note
that $\Lambda$ only depends on the imaging setup and the choice of discretization of the IW.
For simplicity in \eqref{eq:ptransform}, we have fixed the coordinate system so that $x_s=s\, a/N$, $s=0,\dots,N-1$, 
$x_r=r\, a/N$, $r=0,\dots,N-1$, $x'_j=j\,b/K$, $j=1,\dots,K$. 

Clearly,  in the paraxial approximation, imaging is closely related to a discrete Fourier transform. In order to realize this connection, 
fix $s=0$  and view 
$\vect H = \vect H( \tilde{\rho} )$ as a linear map from the IW to the array, so
\begin{equation}
\vect H: \tilde{\rho}  \to  \sum_{j=1}^{K} \tilde{\rho}_j e^{- i \frac{ \Lambda  r j}{N}}, r=0,\dots, N-1.
\label{eq:H}
\end{equation}
The above is a Fourier transform if $N=K$ and $\Lambda= 2\pi$, which holds when the number of discrete points  in the IW is
 $K=K_{opt}=a b (\lambda L)^{-1}=F\,b/a$, where $F$ is the Fresnel number \eqref{eq:fresnelnumber}.
We refer to $K_{opt}$ as the optimal number  of sampling points in the IW because for $N=K=K_{opt}$ \eqref{eq:H}
is an orthogonal transformation. Note that $K_{opt}$ is  the number of points in the IW if  the sampling interval is the well 
known cross-range resolution limit $\lambda L/a$. 

Note that since $\vect H$ is a Hankel matrix, only two illuminations (from both edges of the array) are needed to obtain the full response matrix, 
if the phases can be recorded. 
All the other illuminations are redundant in the one-dimensional setup when the paraxial approximation holds and the scatterers are all at the same range. 

Next, we show that six illuminations are enough to recover $\vect H$ completely from intensity-only measurements. 
Let us assume that the phases cannot be measured and, therefore, only the intensities of the received signals are available for imaging. 
Let
\[
\vect f^{(1)}=\vect e_1, \,\,
\vect f^{(2)}= \vect D^{-1}_{\mbox{\scriptsize source}} \left( \vect e_1+ \vect e_2 \right),\,\, \mbox{and}\,\,
\vect f^{(3)}=\vect D^{-1}_{\mbox{\scriptsize source}} \left( \vect e_1+ i \vect e_2 \right)\, ,
\]
with 
$
\vect e_1=[ 1,0,0\dots,0]^t$ and $\vect e_2=[ 0,1,0\dots,0]^t, 
$
be three illuminations that use the two top transducers of the array. When we use the above illuminations we measure
\begin{equation}\label{b1}
\diag \left( \left(\vect \Pm \vect f^{(1)} \right) \left( \vect \Pm \vect f^{(1)} \right)^* \right)
=\diag \left( \left(\vect H \vect e_{1} \right) \left( \vect H \vect e_{1} \right)^* \right) = \abs{ \vect{h^{(1)}} }^2,
\end{equation}
\begin{equation}\label{b2}
\diag \left( \left(\vect \Pm \vect f^{(2)} \right) \left( \vect \Pm \vect f^{(2)} \right)^* \right)
 = \abs{ \vect{h^{(1)}} +\vect{h^{(2)}}}^2,
 \end{equation}
and
\begin{equation}\label{b3}
 \diag \left( \left(\vect \Pm \vect f^{(3)} \right) \left( \vect \Pm \vect f^{(3)} \right)^* \right)
 = \abs{ \vect{h^{(1)}} + i \vect{h^{(2)}}}^2\, ,
\end{equation}
respectively. In \eqref{b1}-\eqref{b3}, $\vect{h^{(1)}}$ and $\vect{h^{(2)}}$ are the first and second columns of the matrix $\vect H$.
Since $\vect H$ is a Hankel matrix, \eqref{b1},~\eqref{b2}, and~\eqref{b3} are
\[
 \abs{ \vect{h^{(1)}} }^2 = \left[ \abs{\Xi_1}^2, \abs{\Xi_2}^2, \dots, \abs{\Xi_{N}}^2 \right]^t,  \qquad  \qquad \quad
\]
\[
 \abs{ \vect{h^{(1)}} +\vect{h^{(2)}}}^2=\left[ \abs{\Xi_1 +\Xi_2 }^2, \abs{\Xi_2 + \Xi_3 }^2, \dots,  
\abs{\Xi_N + \Xi_{N+1} }^2 \right]^t, 
\]
and
\[
\abs{ \vect{h^{(1)}} + i \vect{h^{(2)}}}^2= \left[  \abs{\Xi_1 + i \Xi_2 }^2, \abs{\Xi_2 + i \Xi_3 }^2, \dots,  
\abs{\Xi_N + i \Xi_{N+1} }^2\right]^t. 
\]
Using the polarization identity
we can, therefore, obtain the  dot products
\[
\langle \Xi_i,  \Xi_{i+1} \rangle, \quad i=1,\dots,N\, ,
\]
from the intensities~\eqref{b1}-\eqref{b3}. Since $\abs{\Xi_i}$ are known for  $i=1,2,\dots, N$,  
 we can retrieve $\Xi_i$, $i=1,2,\dots, N$ up to a common phase. Similarly, we can use the following three different illuminations 
 \[
\vect f^{(4)}=\vect e_N, \vect f^{(5)}= \vect D^{-1}_{\mbox{\scriptsize source}} \left( \vect e_N+ \vect e_{N-1} \right),
\vect f^{(6)}=\vect D^{-1}_{\mbox{\scriptsize source}} \left( \vect e_{N}+ i \vect e_{N-1} \right)\, ,
\]
from the two transducers at the bottom  of the array, to find the remaining $\Xi_i$ for  $i=N+1,2,\dots, 2 N-1$. 
If, for definiteness, we assume that $\Xi_1$ is a real positive number, 
then  $\Xi_i$ $i=1,2,\dots, 2 N-1$ are determined uniquely, and therefore the matrix $\vect H$ is reconstructed. 
Our construction again relies on the assumption that all $\Xi_i \neq 0$, which, as we mentioned ealier~\cite{DS}, is not restrictive. 

The two-dimensional case consisting of a planar array and a planar image is analogous to the one described above. 
For a  two-dimensional setup we will need to make twelve measurements to recover the full-response matrix instead of six as we showed 
in the one-dimensional setup.

\section{Numerical Simulations}
\label{sec:numerics}

\begin{figure}[t]
    \centering
    \begin{tikzpicture}[scale=0.45, transform shape]
       \node at (3.5, 7.1) {\huge{(a)}};
        \draw [<-] (0, -0.1) -- (0, 3.2);
        \node at (0, 3.5) {\Large{$a$\Large}};
        \draw [->] (0, 3.9) -- (0, 7.1);
        \fill[fill=black] (1.5, 0) circle [radius=0.1];
        \fill[fill=black] (1.5, 1) circle [radius=0.1];
        \fill[fill=black] (1.5, 2) circle [radius=0.1];
        \node at (1, 2) {\Large{$\vect x_s$}};
        \fill[fill=black] (1.5, 3) circle [radius=0.1];
        \fill[fill=black] (1.5, 4) circle [radius=0.1];
        \fill[fill=black] (1.5, 5) circle [radius=0.1];
        \fill[fill=black] (1.5, 6) circle [radius=0.1];
        \node at (1, 6) {\Large{$\vect x_r$}};
        \fill[fill=black] (1.5, 7) circle [radius=0.1];
        \draw [<-] (1.7, -0.3) -- (5.2, -0.3);
        \node at (5.6, -0.3) {$L$};
        \draw [->] (6.0, -0.3) -- (9.3, -0.3);
        \fill[fill=black] (8.9, 3.7) circle [radius=0.1];
        \node at (9.6, 3.1) {\Large{$\vect y_{n_1}$}};
        \fill[fill=black] (9.3, 3.5) circle [radius=0.1];
        \node at (8.5, 3.5) {\Large{$\vect y_{n_2}$}};
        \fill[fill=black] (9.5, 3.9) circle [radius=0.1];
        \node at (9.4, 4.3) {\Large{$\vect y_{n_3}$}};
        \node at (10, 2.1) {\Large{$IW$}};
        \draw [-] (8, 2.5) -- (8, 4.8);
        \draw [-] (8, 4.8) -- (10, 4.8);        
        \draw [-] (10, 2.5) -- (10, 4.8); 
        \draw [-] (8, 2.5) -- (10, 2.5);    
        \draw[-latex,line width=0.6pt] (2, 2) -- (8.5, 3) node[pos=0.5, sloped, above]{\Large$\wG(\vect x_s, \vect y)$};
        \draw (1.7, 2) arc[radius=0.2, start angle=0, end angle=30];
        \draw (1.7, 2) arc[radius=0.2, start angle=0, end angle=-30];
        \draw (2, 2) arc[radius=0.5, start angle=0, end angle=30];
        \draw (2, 2) arc[radius=0.5, start angle=0, end angle=-30];
        \draw (2.3, 2) arc[radius=0.9, start angle=0, end angle=30];
        \draw (2.3, 2) arc[radius=0.9, start angle=0, end angle=-30];
        \draw (2.6, 2) arc[radius=1.2, start angle=0, end angle=30];
        \draw (2.6, 2) arc[radius=1.2, start angle=0, end angle=-30];
        \draw (2.9, 2) arc[radius=1.5, start angle=0, end angle=30];
        \draw (2.9, 2) arc[radius=1.5, start angle=0, end angle=-30];
    \end{tikzpicture}
\hspace{1.3cm}
    \begin{tikzpicture}[scale=0.45, transform shape]
       \node at (3.5, 7.1) {\huge{(b)}};
        \fill[fill=black] (1.5, 0) circle [radius=0.1];
        \fill[fill=black] (1.5, 1) circle [radius=0.1];
        \fill[fill=black] (1.5, 2) circle [radius=0.1];
        \node at (1, 2) {\Large{$\vect x_s$}};
        \fill[fill=black] (1.5, 3) circle [radius=0.1];
        \fill[fill=black] (1.5, 4) circle [radius=0.1];
        \fill[fill=black] (1.5, 5) circle [radius=0.1];
        \fill[fill=black] (1.5, 6) circle [radius=0.1];
        \node at (1, 6) {\Large{$\vect x_r$}};
        \fill[fill=black] (1.5, 7) circle [radius=0.1];
        \draw [<-] (1.7, -0.3) -- (5.2, -0.3);
        \node at (5.6, -0.3) {$L$};
        \draw [->] (6.0, -0.3) -- (9.0, -0.3);
        \fill[fill=black] (9, 3.3) circle [radius=0.1];
        \node at (9.6, 3.3) {\Large{$\vect y_{n_1}$}};
        \fill[fill=black] (9, 3.7) circle [radius=0.1];
        \node at (9.6, 3.7) {\Large{$\vect y_{n_2}$}};
        \fill[fill=black] (9, 4.9) circle [radius=0.1];
        \node at (9.6, 4.9) {\Large{$\vect y_{n_3}$}};
        \node at (10, 2.1) {\Large{$IW$}};      
        \draw [-] (9, 2.2) -- (9, 5.5); 
        \draw [-] (8.8, 2.2) -- (9.2, 2.2); 
        \draw [-] (8.8, 5.5) -- (9.2, 5.5); 
        \draw[-latex,line width=0.6pt] (2, 2) -- (8.5, 3) node[pos=0.5, sloped, above]{\Large$\wG(\vect x_s, \vect y)$};
        \draw (1.7, 2) arc[radius=0.2, start angle=0, end angle=30];
        \draw (1.7, 2) arc[radius=0.2, start angle=0, end angle=-30];
        \draw (2, 2) arc[radius=0.5, start angle=0, end angle=30];
        \draw (2, 2) arc[radius=0.5, start angle=0, end angle=-30];
        \draw (2.3, 2) arc[radius=0.9, start angle=0, end angle=30];
        \draw (2.3, 2) arc[radius=0.9, start angle=0, end angle=-30];
        \draw (2.6, 2) arc[radius=1.2, start angle=0, end angle=30];
        \draw (2.6, 2) arc[radius=1.2, start angle=0, end angle=-30];
        \draw (2.9, 2) arc[radius=1.5, start angle=0, end angle=30];
        \draw (2.9, 2) arc[radius=1.5, start angle=0, end angle=-30];
    \end{tikzpicture}
    \caption{Schematic. A linear array probes a homogeneous medium sending a spherical wave from $\vect x_s$. Wave propagation is described by the Green's function $\wG(\vect x_s, \vect y)$.
    In figure (a) there are three point-like scatterers at positions $\vect y_{n_j}$, $j=1,2,3$, whose ranges are different and unknown; the IW is a square.
    In figure (b) the point-like scatterers are on a line, and their range is known; the IW is a segment.}
    \label{fig:schematic}
\end{figure}
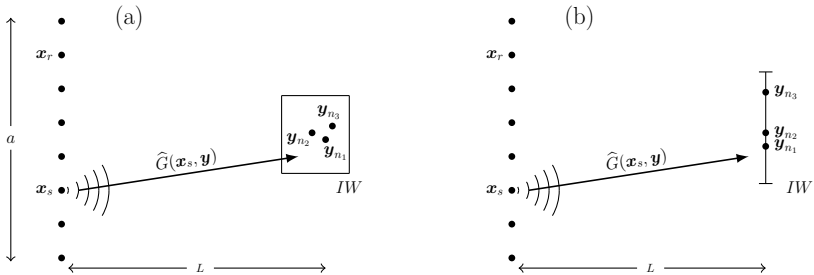

In this section we present numerical simulations that illustrate the performance of the proposed approach, where noise is also included. 
For simplicity of graphical representation we consider a linear array on the $x$ axis. The scatterers lie either on the plane $(x,0,z)$ as in Fig. \ref{fig:schematic} (a), in which case they are at different ranges from the
array, or they lie on the line $(x,0,L)$ as in Fig. \ref{fig:schematic} (b), in which case they are located at the same range $L$ that we consider to be known. The coordinate system has origin at the center of the array and range axis $z$ orthogonal to it. The relevant cross-range coordinate is, therefore, $x$. All the results presented here can be extended to the case in which one considers a planar square array and scatterers with cross-range coordinates $(x,y)$. 

The linear array has aperture size $a=2000 \lambda$, and $N=101$ transducers located at positions $(x_i,0,0)$, $i=1,\dots, N$, that are $20$ wavelengths apart. The scatterers, whose complex reflectivities are set randomly, are placed within a IW whose distance to the array vary from one numerical simulation to another.
The IW is discretized using a uniform lattice with mesh size according to the distance $L$ to the array. We choose
the mesh size to be $h_z=\lambda L^2/a^2$ in the range direction, and $h_x=\lambda L/a$ in the cross-range direction. 
These  are  the well known resolutions limits in array imaging. Therefore, the smaller the distance between the array and the IW,  the smaller the resolution is. To keep the number of grid points independent of the distance to the array we also change the size of the IW accordingly.

In all our numerical simulations  the noise is modeled as additive. More specifically, the
noise at the i-th receiver is modeled by adding a  random variable $\zeta_i$
uniformly distributed  on $[(1-\eps) \beta_i, (1+\eps) \beta_i]$, 
where $\beta_i$ is the noiseless intensity
received on the i-th receiver (see~\eqref{dataI}), and $\eps\in (0,1)$ is a parameter that measures the noise strength. The noise added to the data gathered at different receivers when different illuminations are used is independent.

\begin{figure}[t]
\centering
\begin{tabular}{ccc}
\includegraphics[scale=0.24]{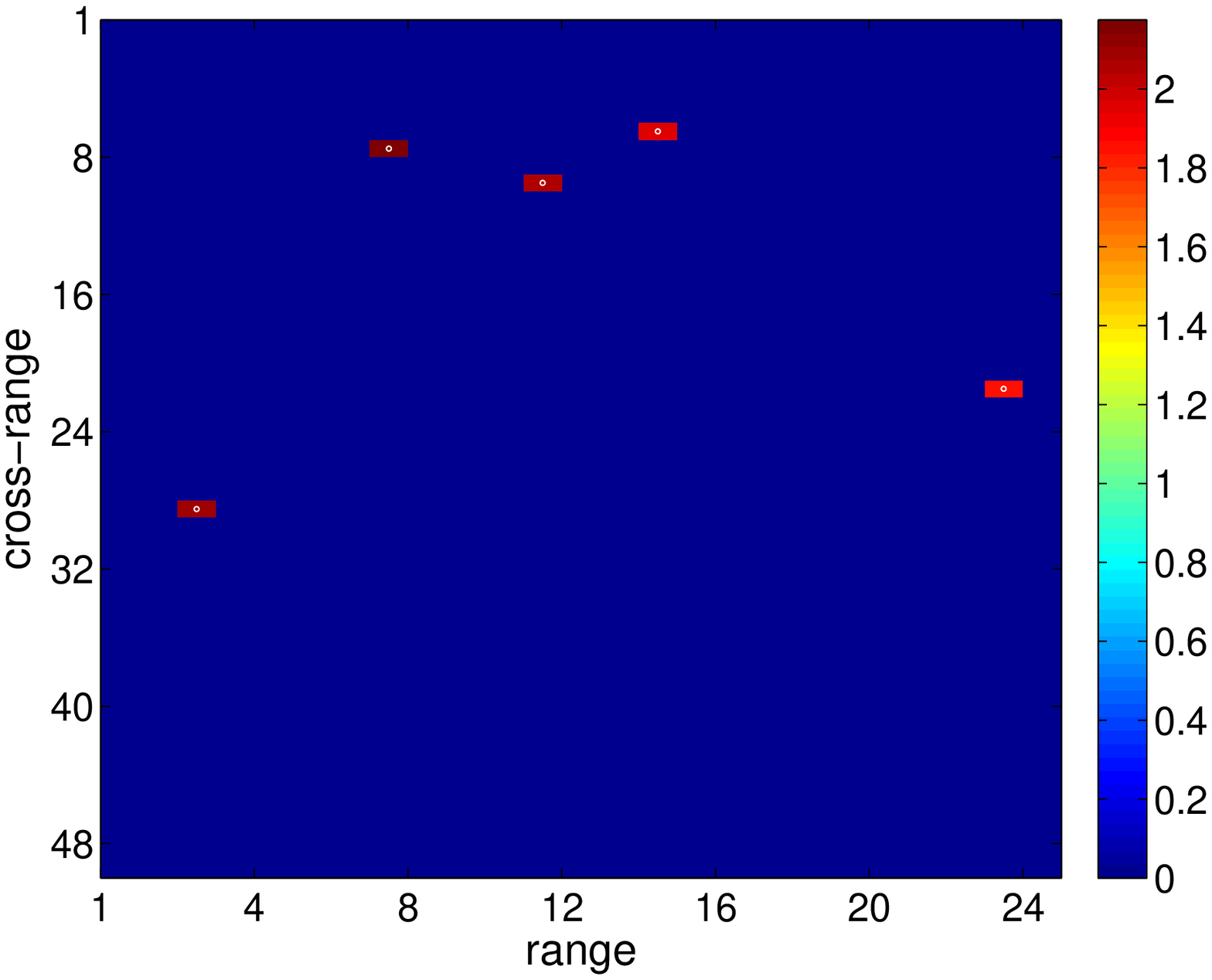} &
\includegraphics[scale=0.24]{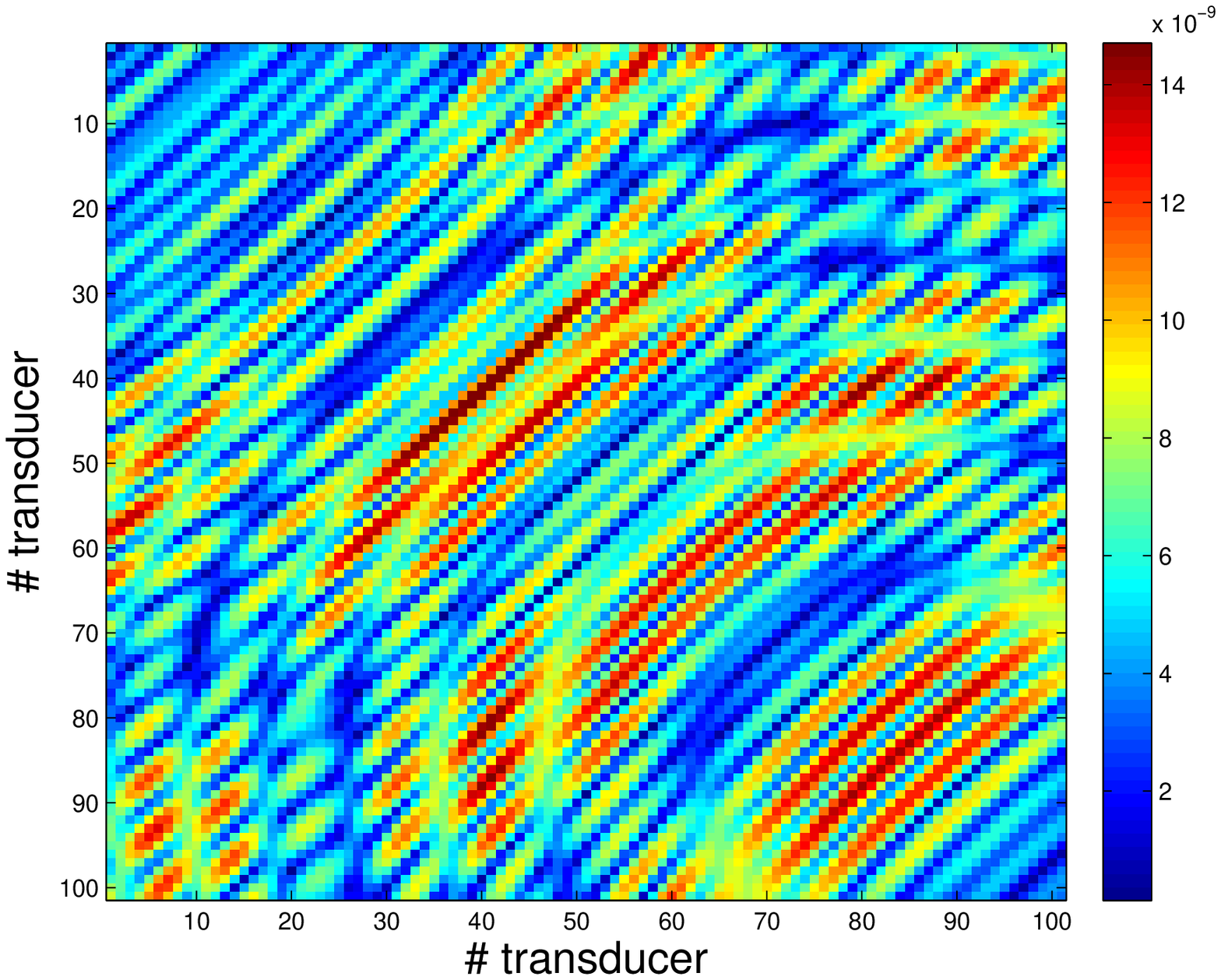} &
\includegraphics[scale=0.24]{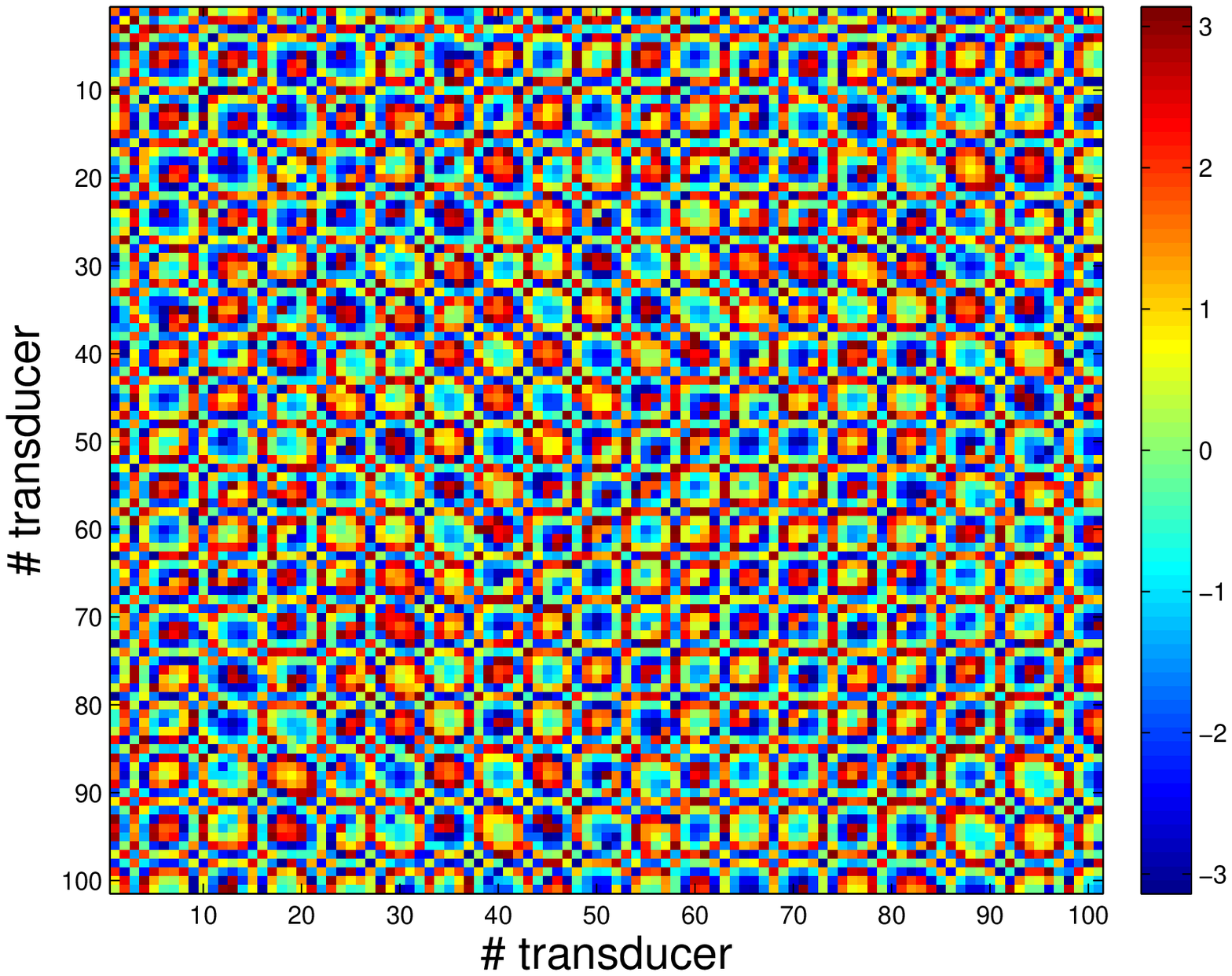} \\
\includegraphics[scale=0.24]{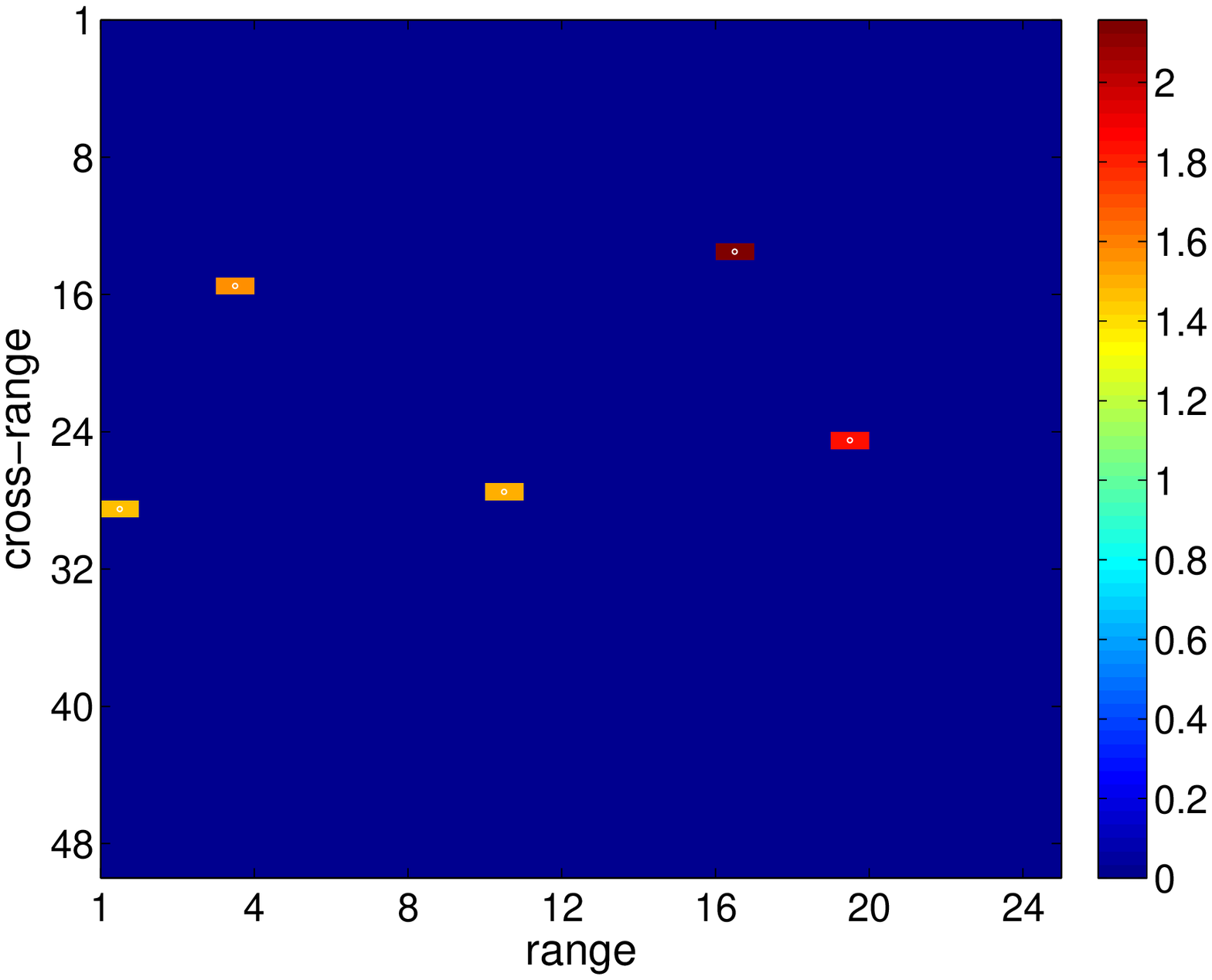} &
\includegraphics[scale=0.24]{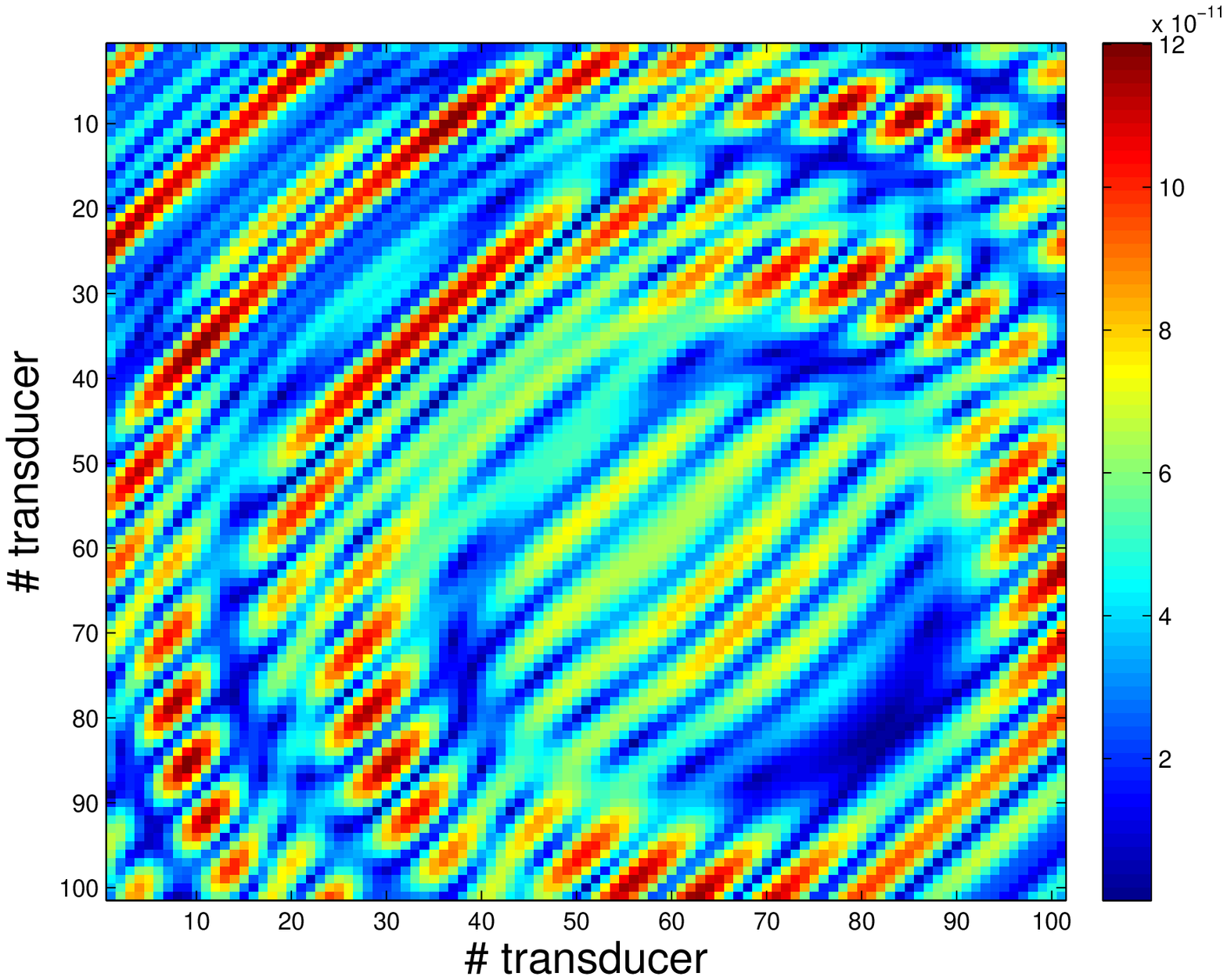} &
\includegraphics[scale=0.24]{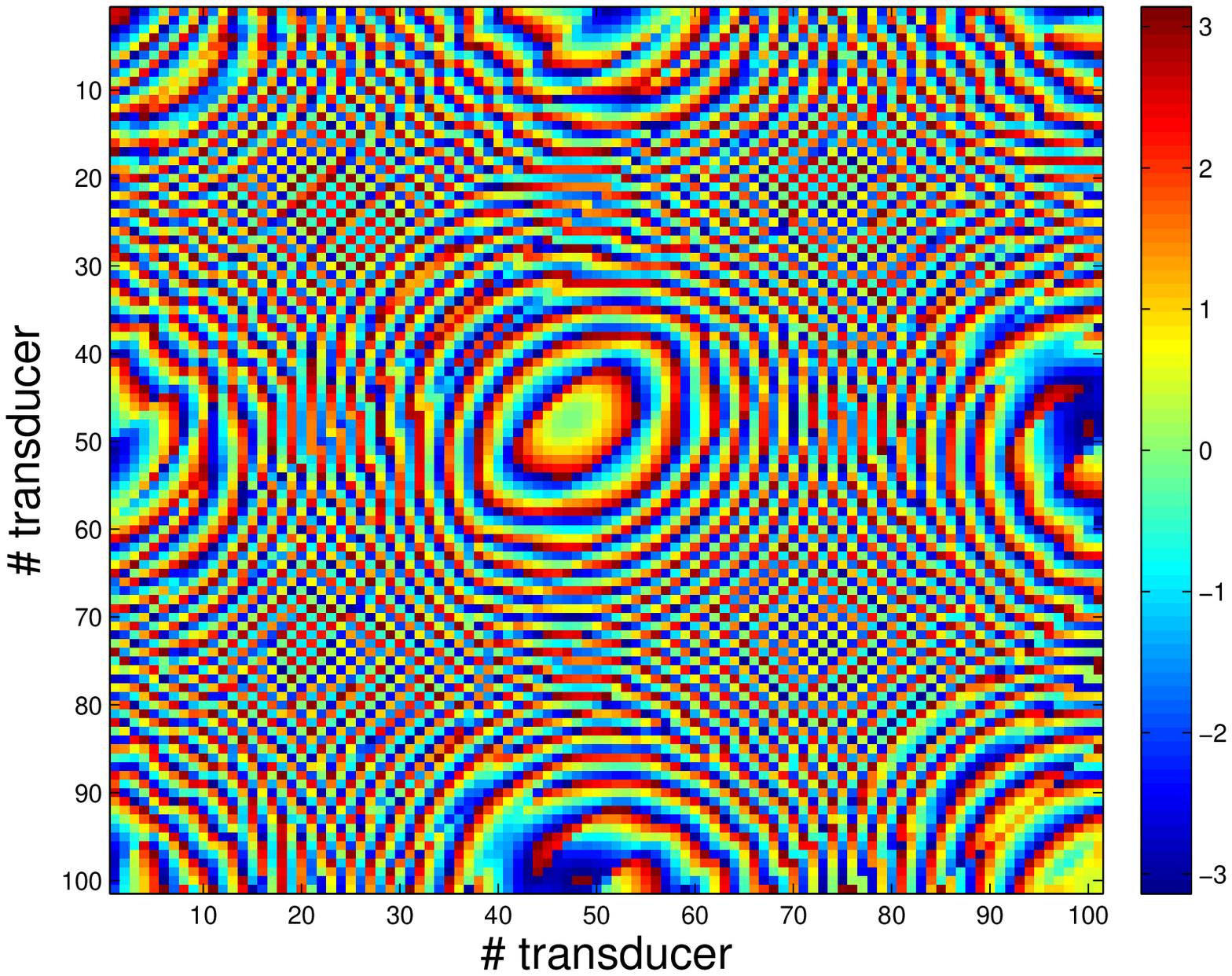} \\
\includegraphics[scale=0.24]{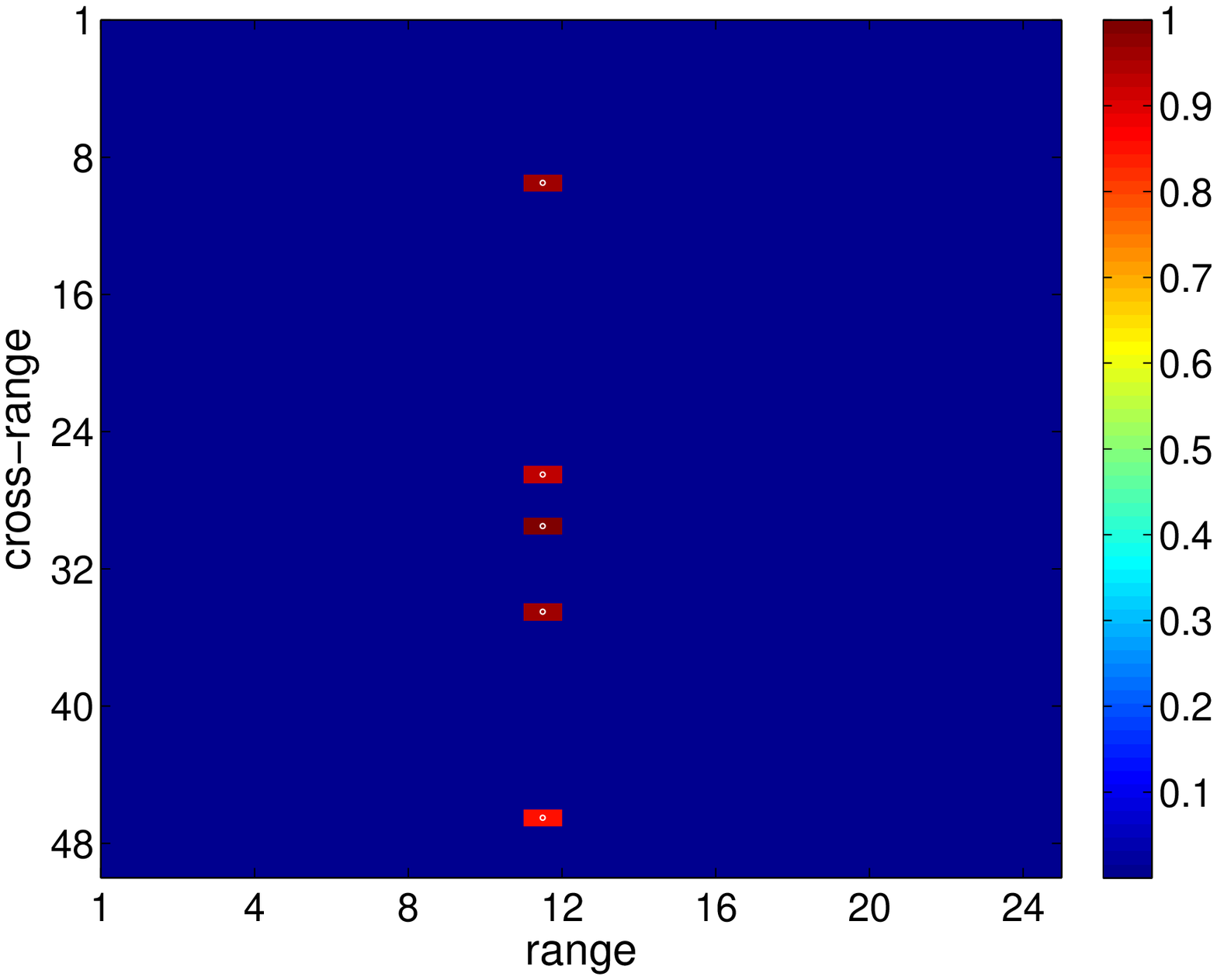} &
\includegraphics[scale=0.24]{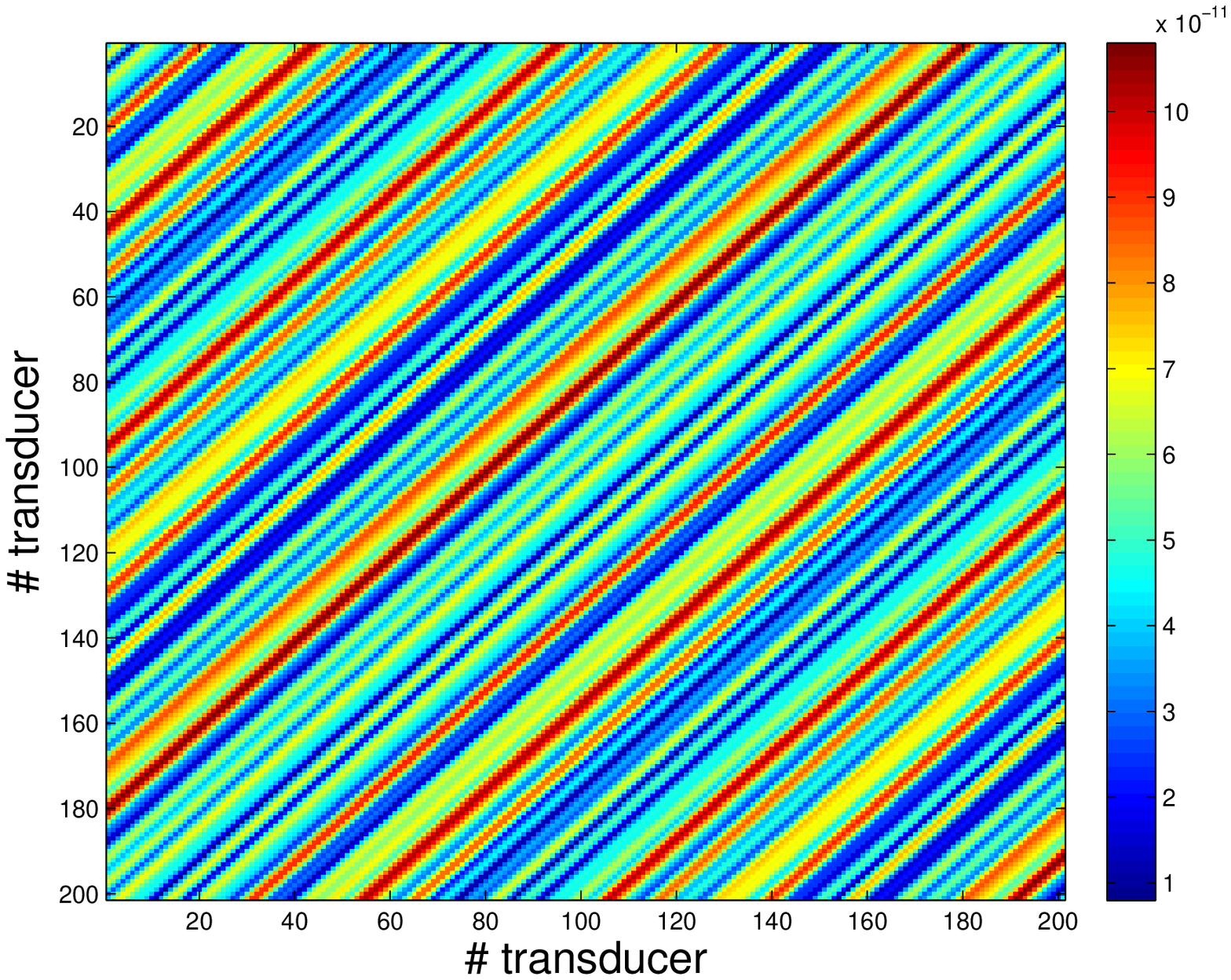} &
\includegraphics[scale=0.24]{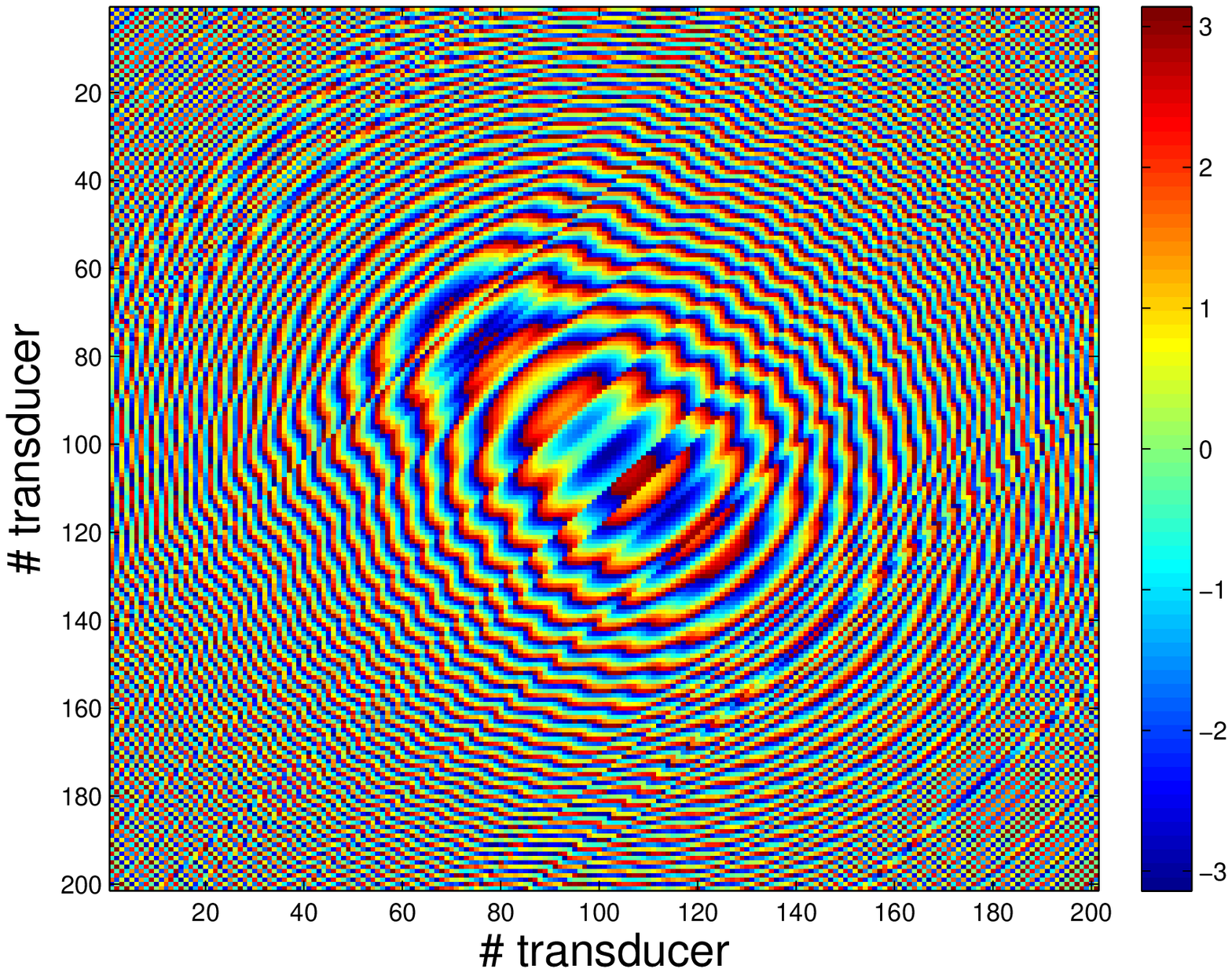} 
\end{tabular}
\caption{Reference images (left column), and amplitudes (middle column) and phases (right column) of the response matrix $P$. 
The distances from the IWs to the array are $L=2000$ (top row) and $L=20000$ (middle and bottom row). In the left column, the range axes are in units of 
$h_z=\lambda L^2/a^2$ and the cross-range axes in units of $h_x=\lambda L/a$.
The length of the array is $2000\lambda$.
\label{fig:examples}
}
\end{figure}

The left column of figure \ref{fig:examples} shows three examples of sets of scatterers at distances $L=2000\lambda$ (top row) and $L=20000\lambda$ 
(middle and bottom row). 
The amplitudes and phases of the full response matrices $\vect P$ are shown in the middle and right columns, respectively. All the information needed 
for imaging is encoded in these figures. However, in our case, phases are not available and, hence,  the array imaging problem is 
to detemine the configuration of the scatterers using only the amplitudes of the entries of $\vect P$, i.e., the middle column. 
For very large $L$ (middle row), the paraxial approximation holds, and the $\vect P$ matrix has a much simpler structure. If, in addition, 
all the scatterers are at the same range 
(bottom row), then the amplitude response displays a Hankel structure. The phase response does not display a Hankel structure here, because the 
geometrical factors~\eqref{eq:factor} (that are known) have not been removed.

\begin{figure}[!htb]
\centering
\includegraphics[scale=0.25]{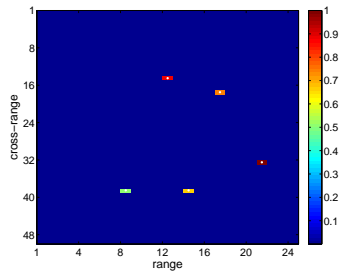} 
\caption{Reference image used in Fig. \ref{fig:range}. The range axis is in units of $h_z=\lambda L^2/a^2$ and the cross-range axis in units of $h_x=\lambda L/a$.
\label{fig:reference}}
\end{figure}

\begin{figure}[!htb]
\centering
\begin{tabular}{ccc}
\includegraphics[scale=0.24]{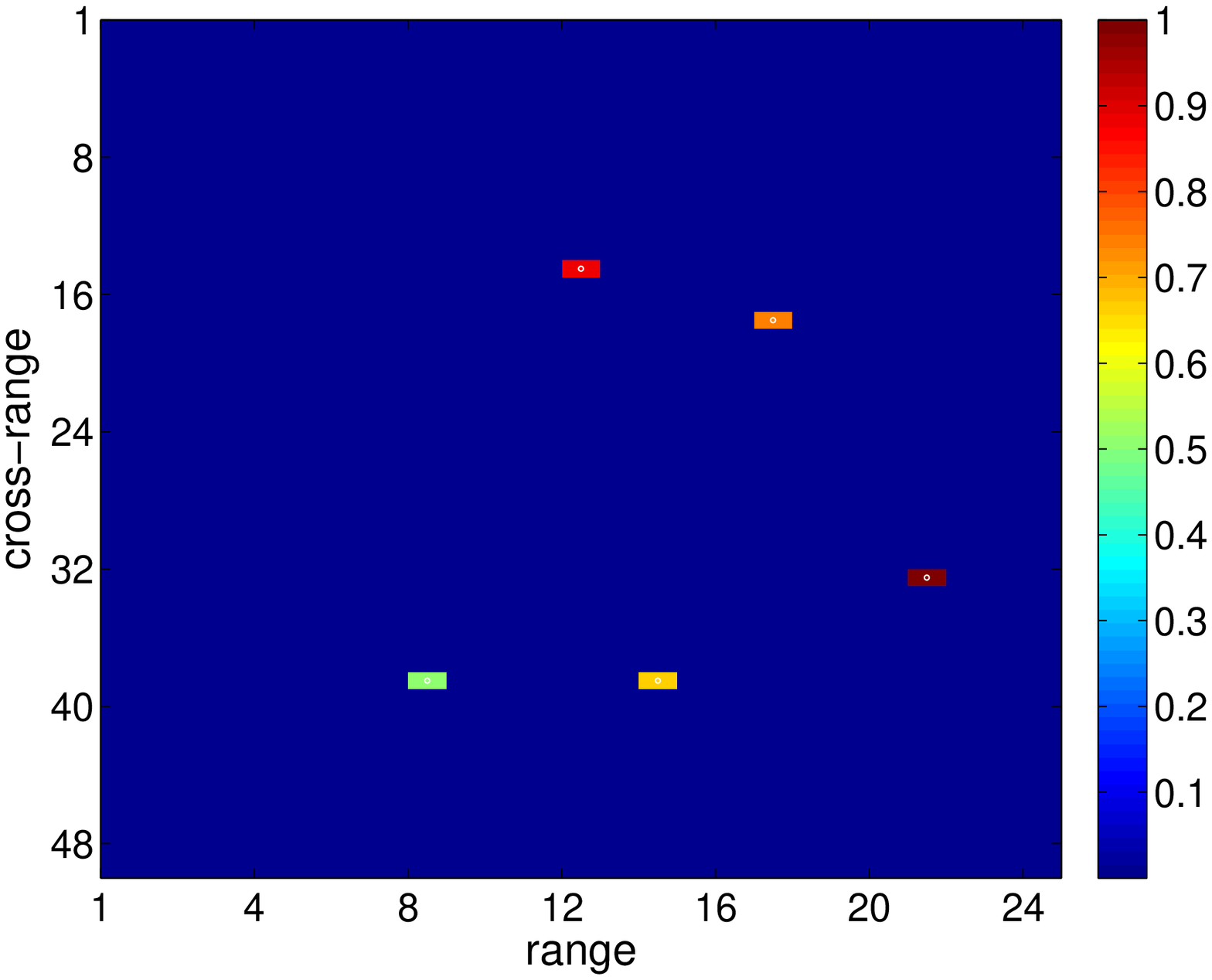} &
\includegraphics[scale=0.24]{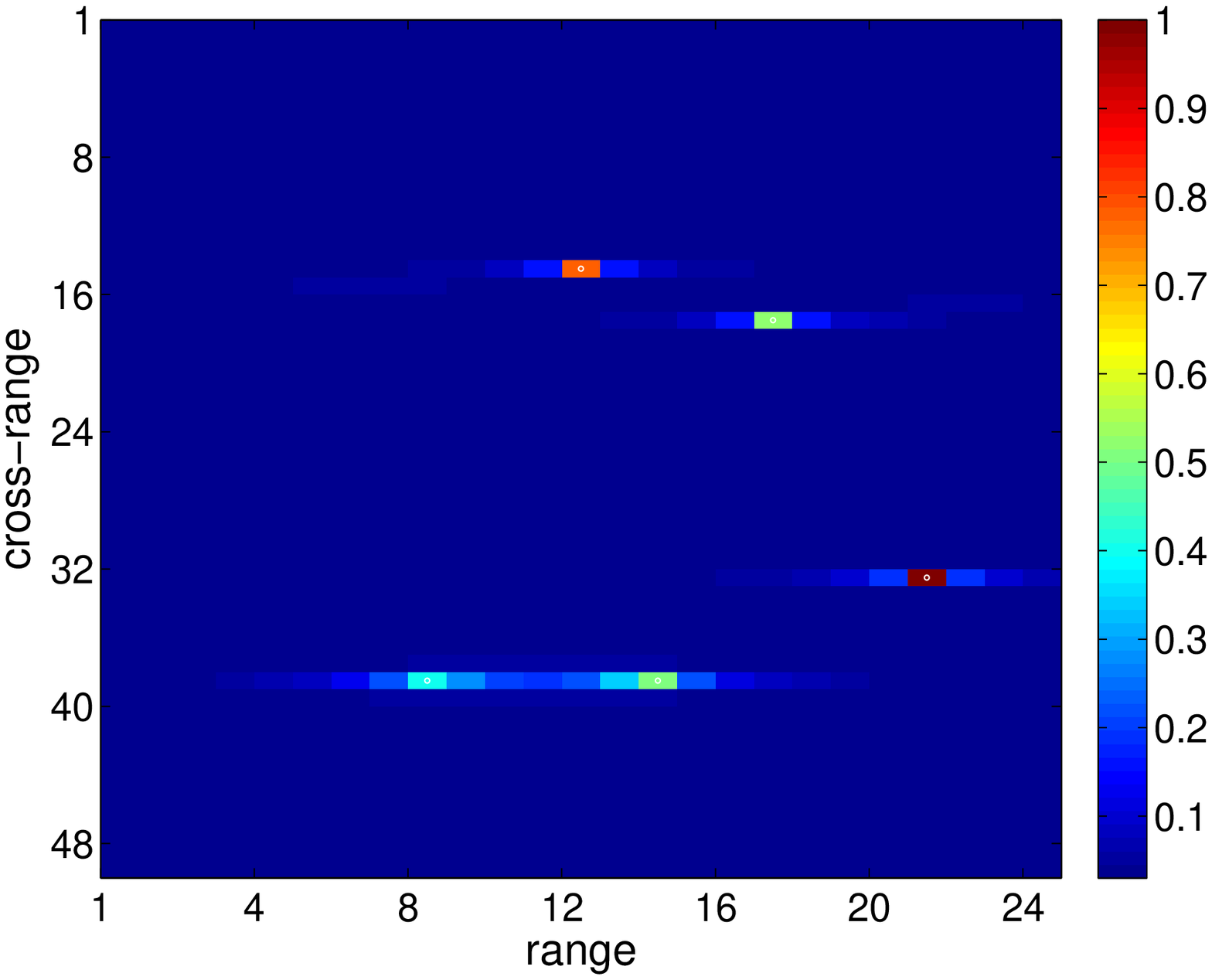} &
\includegraphics[scale=0.24]{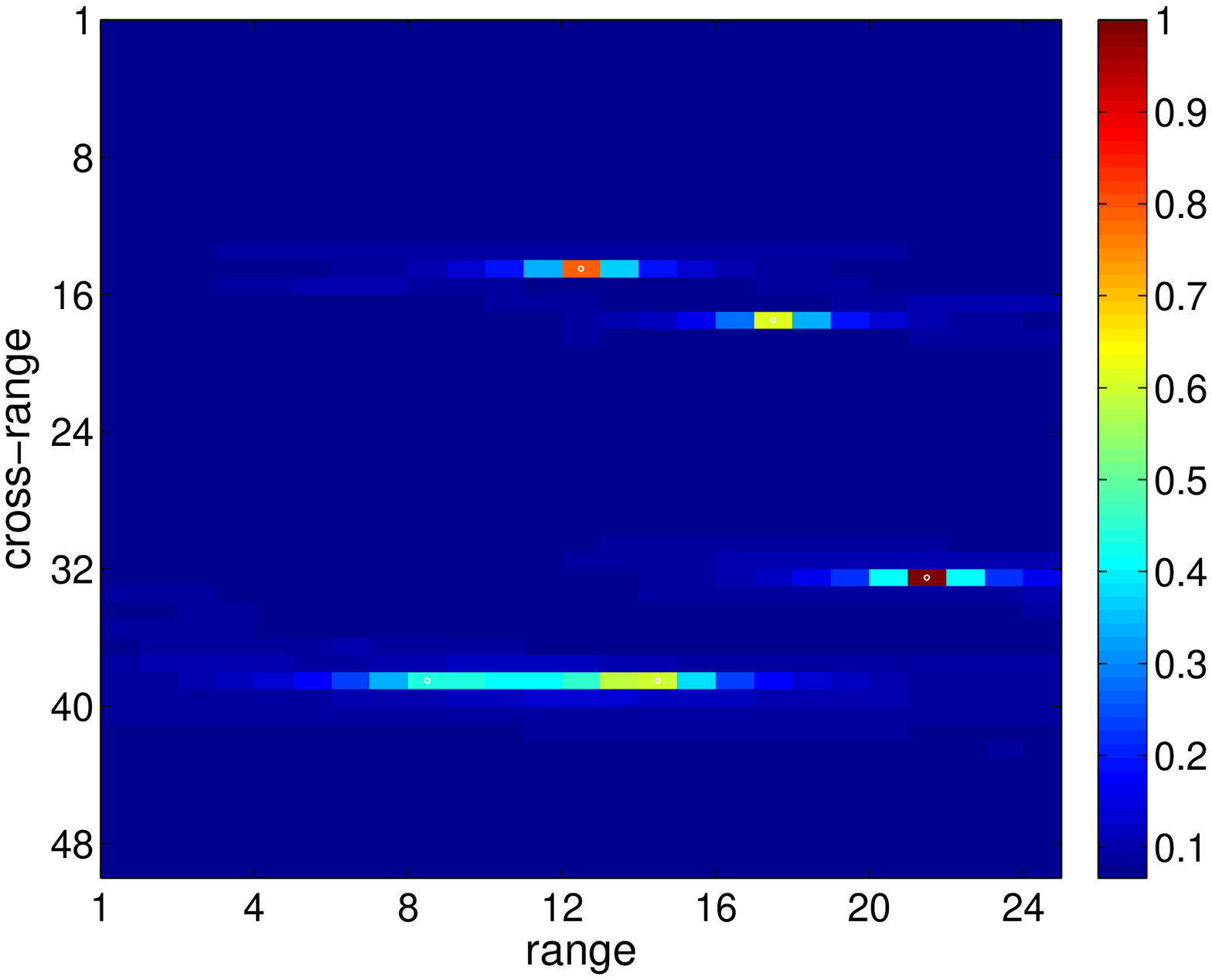} \\
\includegraphics[scale=0.24]{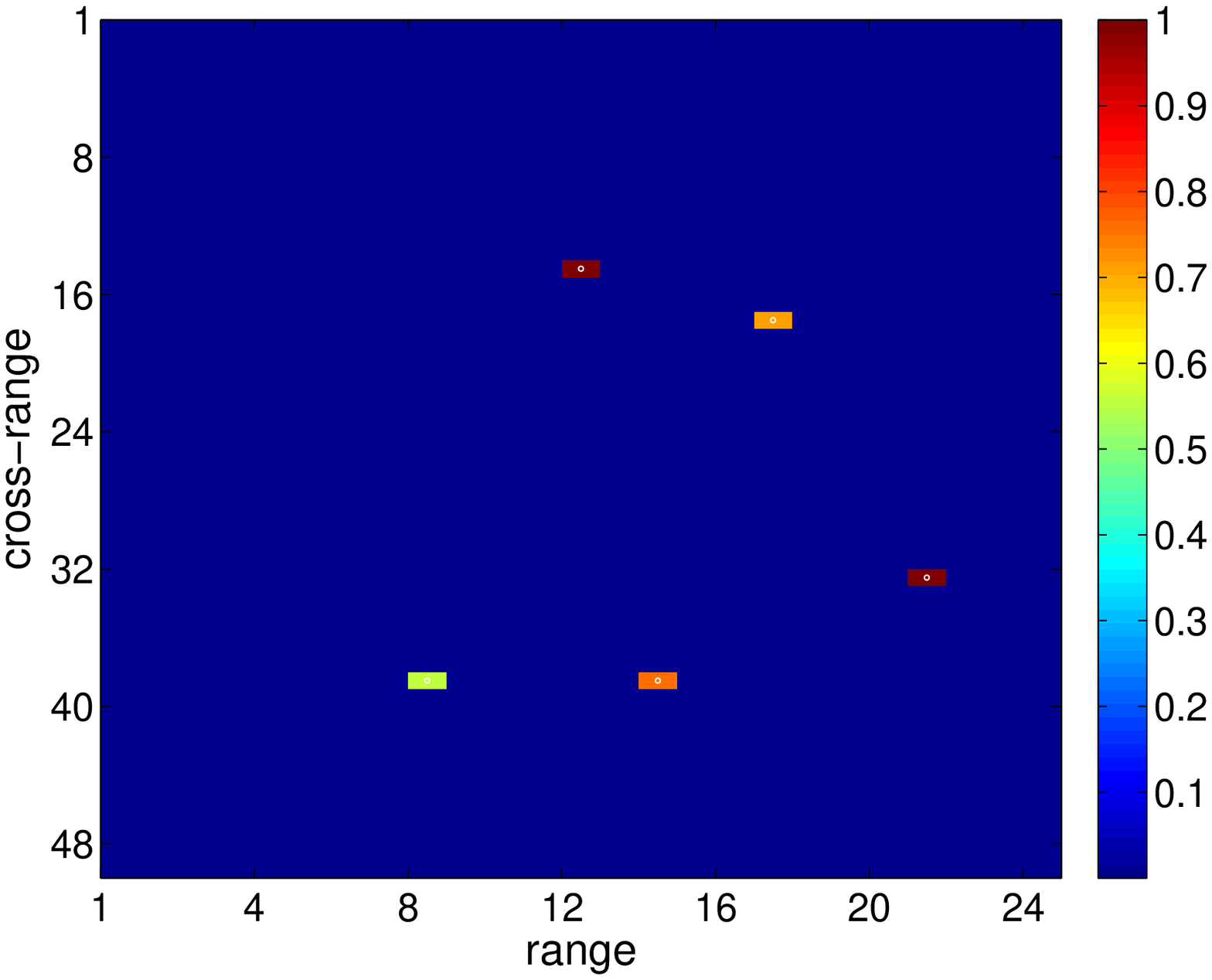} &
\includegraphics[scale=0.24]{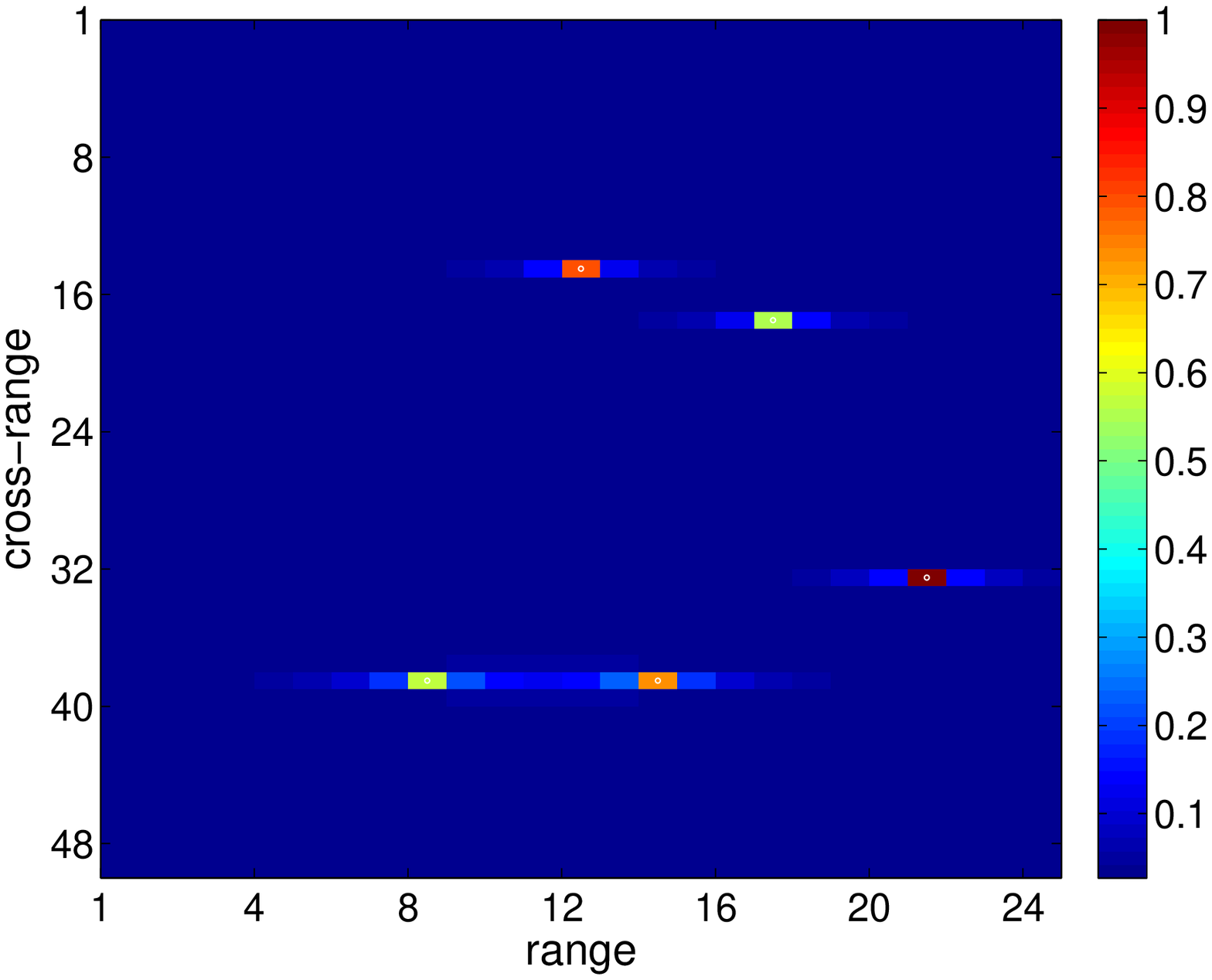} &
\includegraphics[scale=0.24]{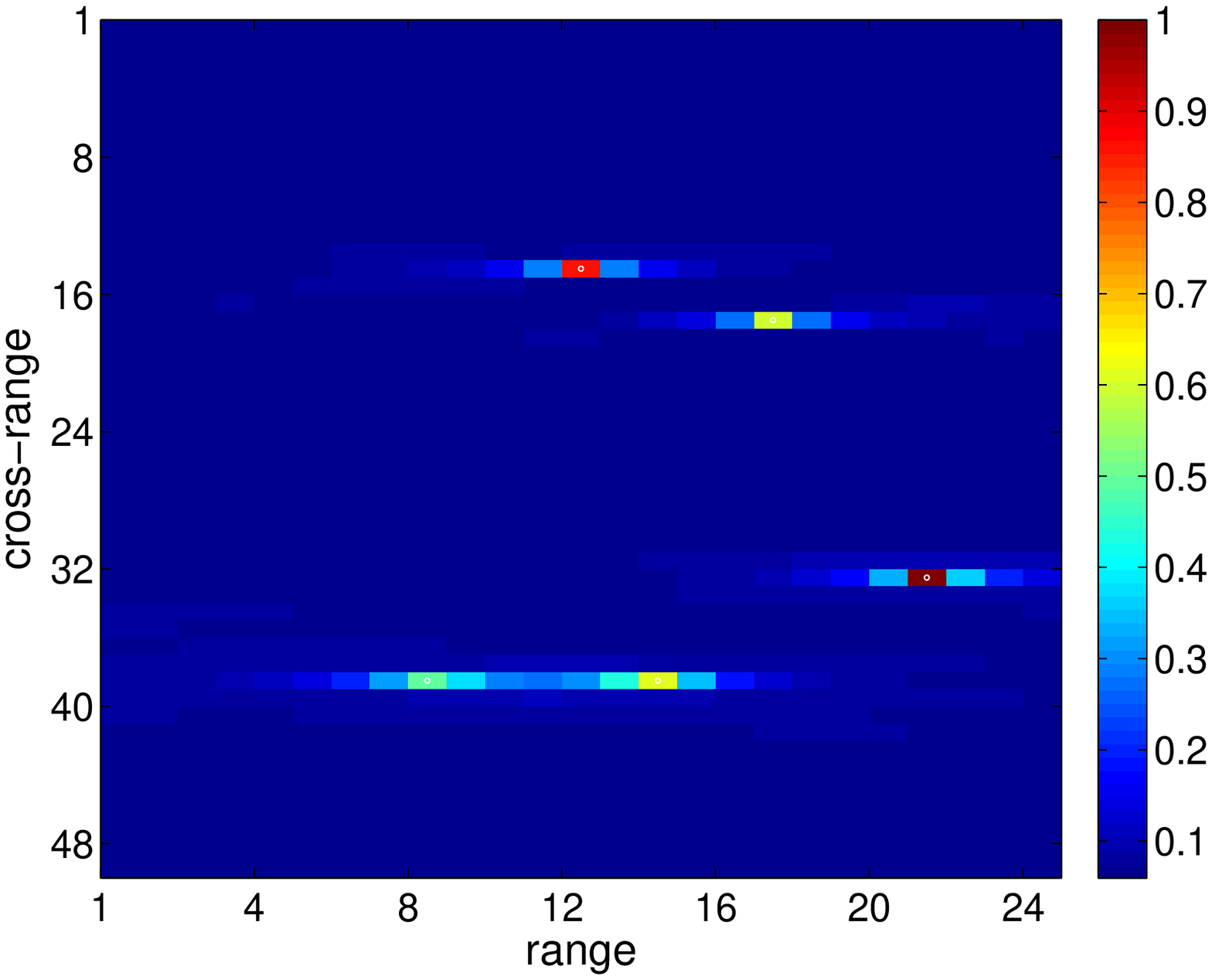} \\
\includegraphics[scale=0.24]{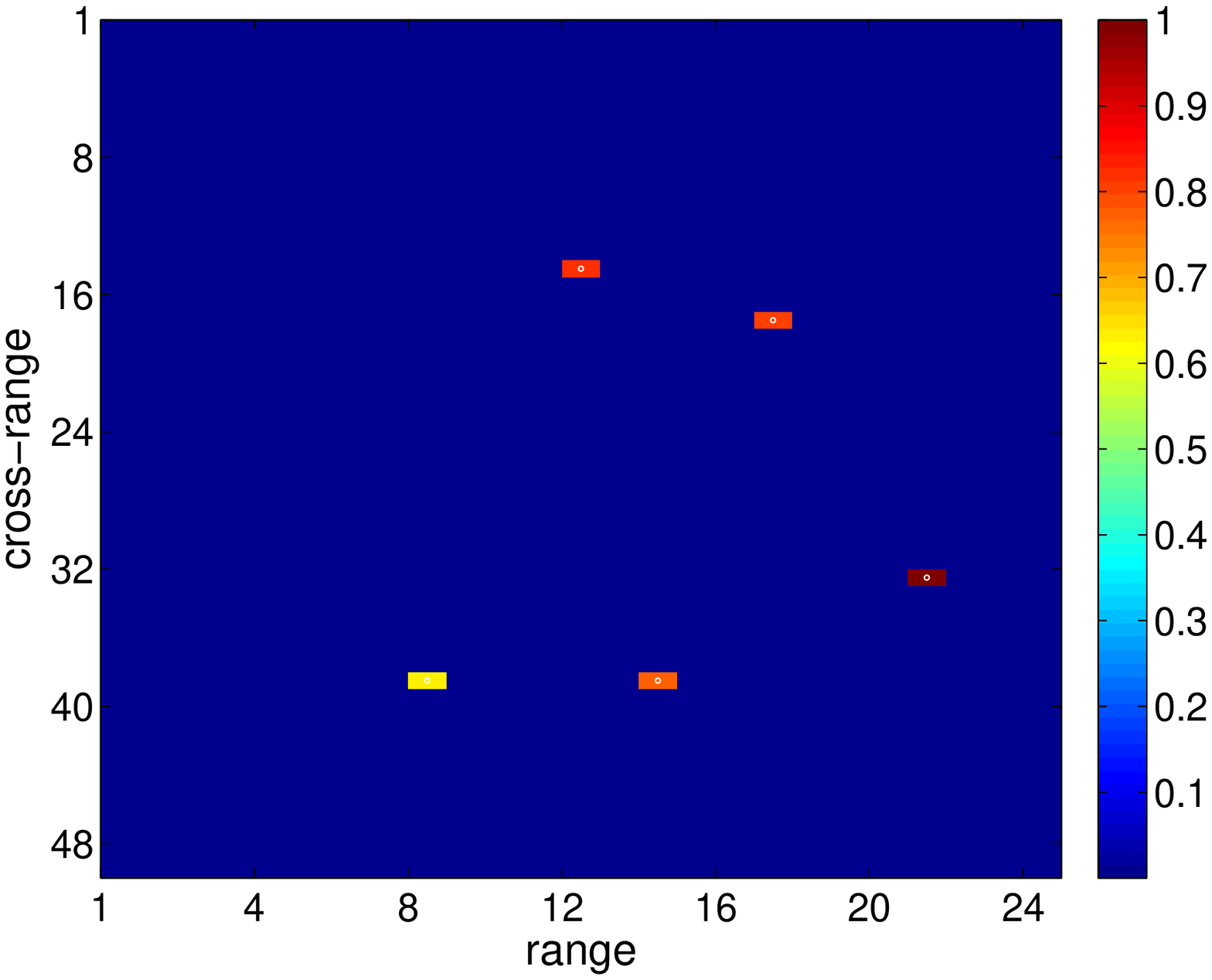} &
\includegraphics[scale=0.24]{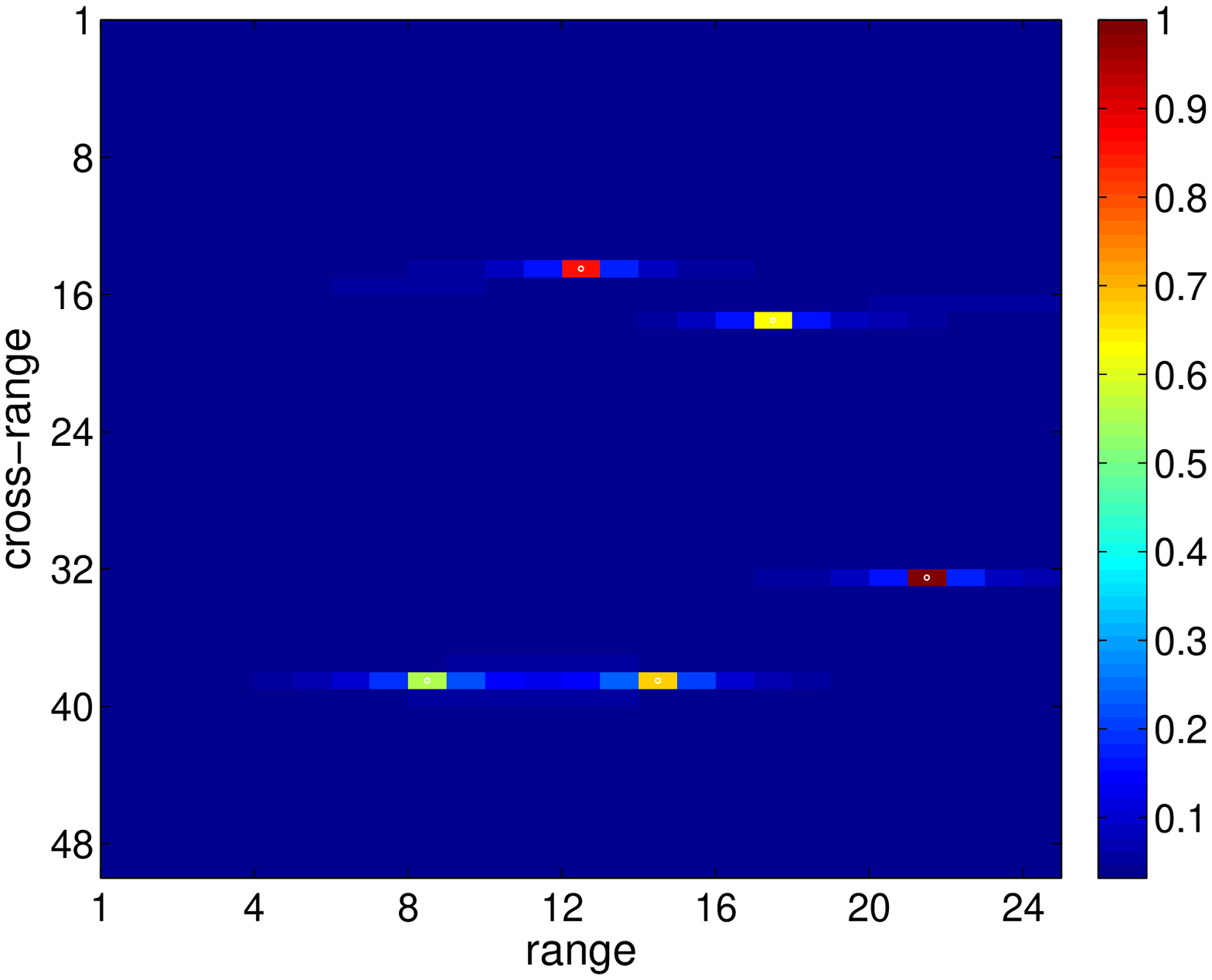} &
\includegraphics[scale=0.24]{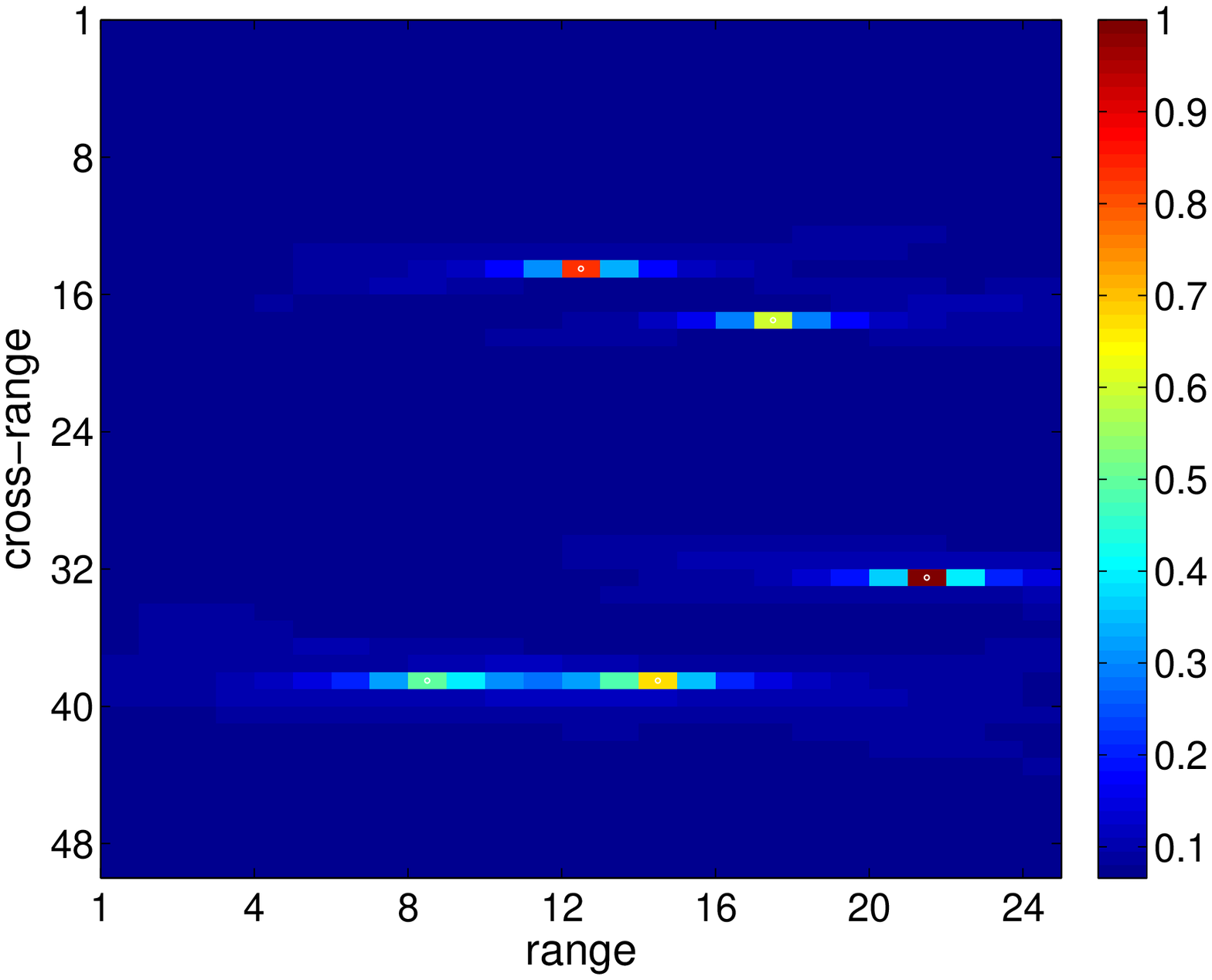} 
\end{tabular}
\caption{Images obtained with MUSIC from the full recovered matrix $M$ for different distances between the IW and the array.
In the top, middle and bottom rows the distances between the IW and the array are $L=2000\lambda$, $L=5000\lambda$ and 
$L=10000 \lambda$, respectively.
The left, middle and right columns contain $0\% $,  $10\% $, $20\% $ of additive noise in the data, respectively. 
Length of the array $a=2000 \lambda$, and number of transducers $Na=101$. The range axes are in units of $h_z=\lambda L^2/a^2$ and the cross-range axes 
in units of $h_x=\lambda L/a$.
\label{fig:range}
}
\end{figure}

Figure \ref{fig:reference} displays the reference IW  used for the numerical simulations  in
Figure \ref{fig:range}. The scatterers are located at different ranges as in Fig. \ref{fig:schematic} (a). 
The range axes is in units of $h_z=\lambda L^2/a^2$, and the cross-range axes in units of $h_x=\lambda L/a$. 
This means that the resolution of the images increases as we move the IW  forward. In Figure~\ref{fig:range}, the distance $L$ between IW and 
the array is $L=2000\lambda$ in the top row, $L=5000\lambda$ in the middle row, and $L=10000\lambda$ in the bottom row.
The images are formed  using the MUSIC (MUltiple
SIgnal Classification) imaging function \cite{Schmidt86}, reviewed in Appendix~\ref{Be}, for the 
full time reversal matrix $\vect M$ obtained with the illumination strategy described in subsection \ref{sec:illum_general}. 
The left column of Figure~\ref{fig:range} shows the results for noise free data. We observe that by using MUSIC for $\vect M$ we get the exact  
locations of the scatterers regardless of the distance between the array and the IW. As expected, when the data is corrupted by $10\%$ (middle column) 
and $20\%$ (right column) of additive noise the images become blurred and the resolution is compromised.
We note, however, that 
all the scatterers are still very well resolved in cross-range.

\begin{figure}[!htb]
\centering
\begin{tabular}{ccc}
\includegraphics[scale=0.24]{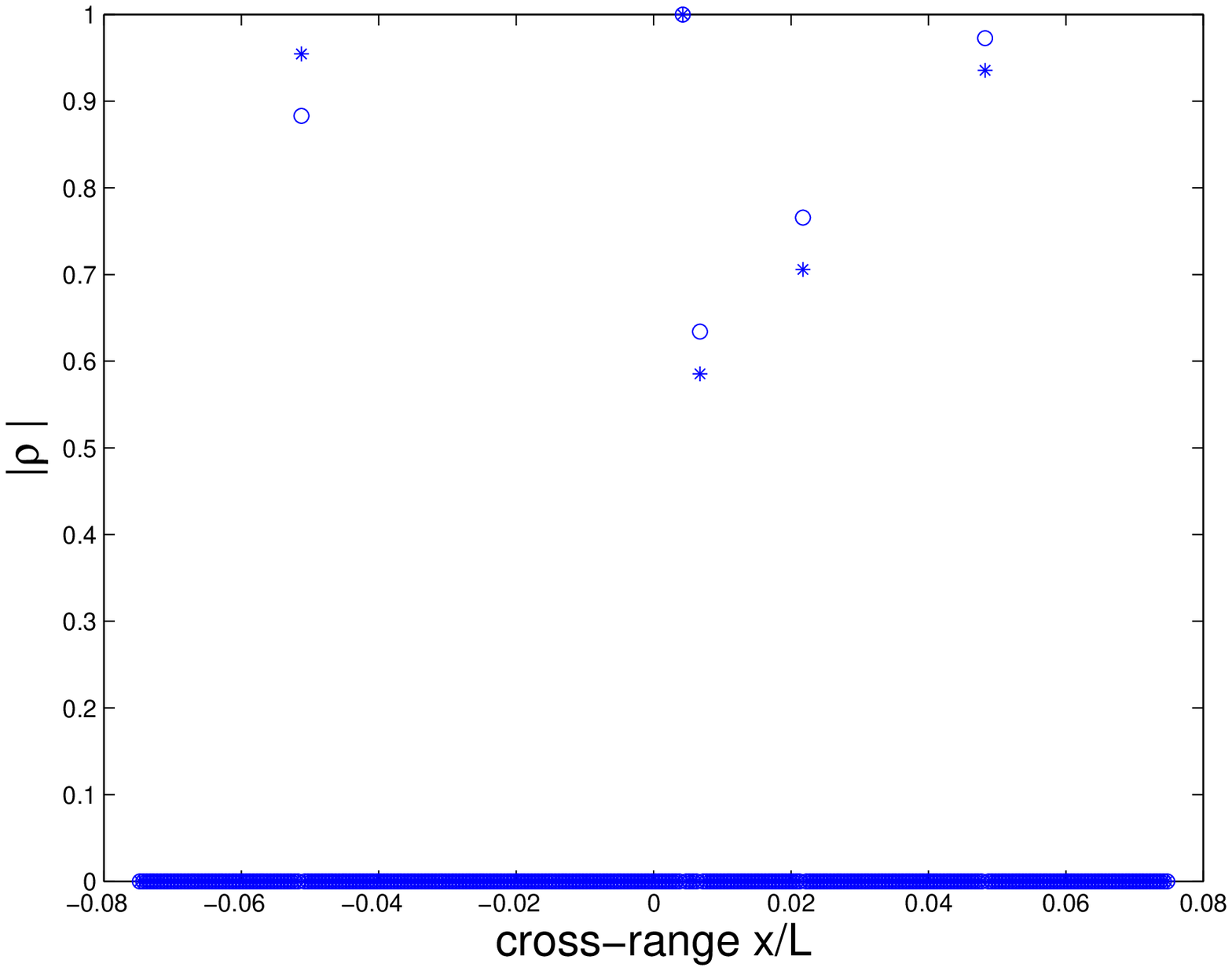} &
\includegraphics[scale=0.24]{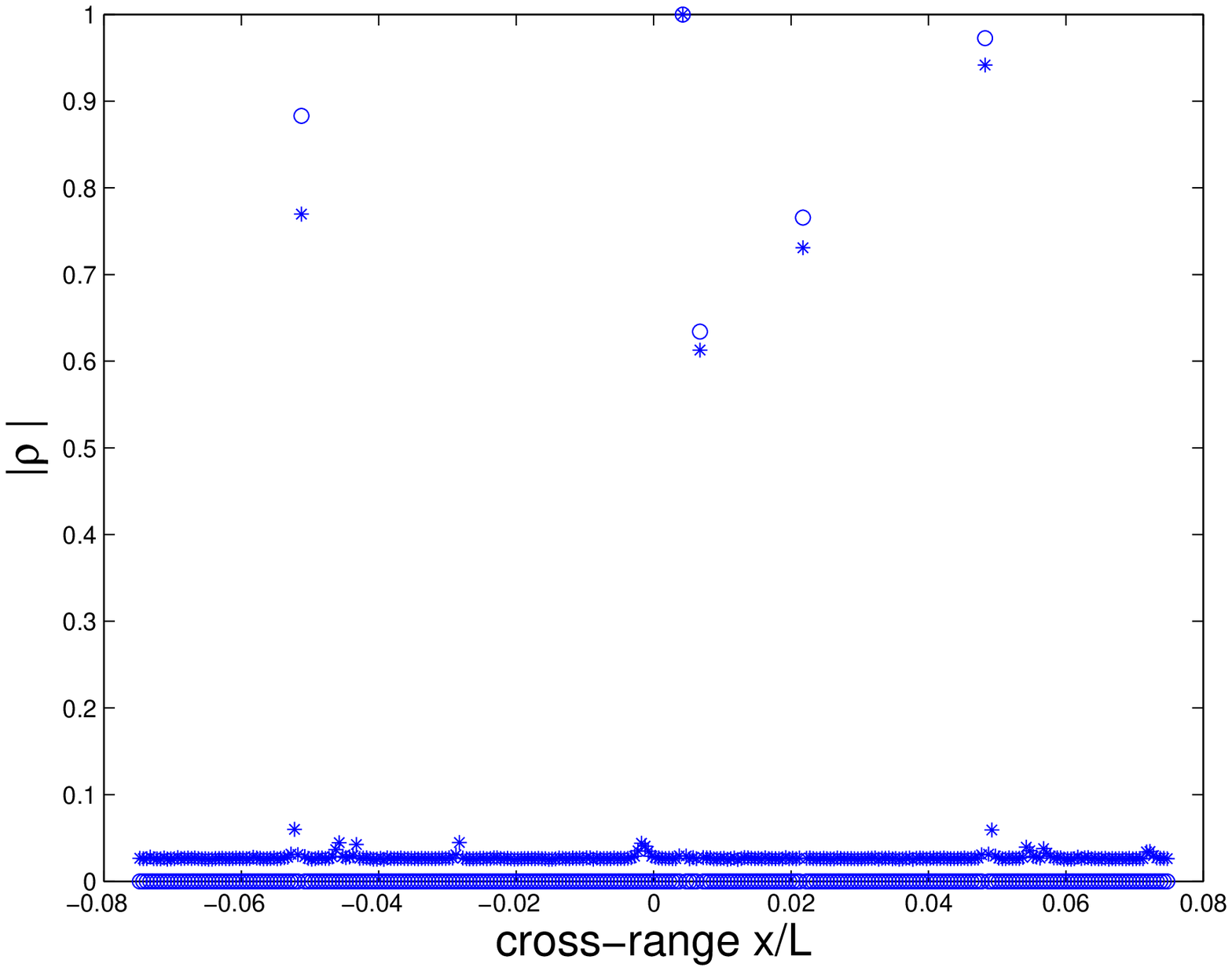} &
\includegraphics[scale=0.24]{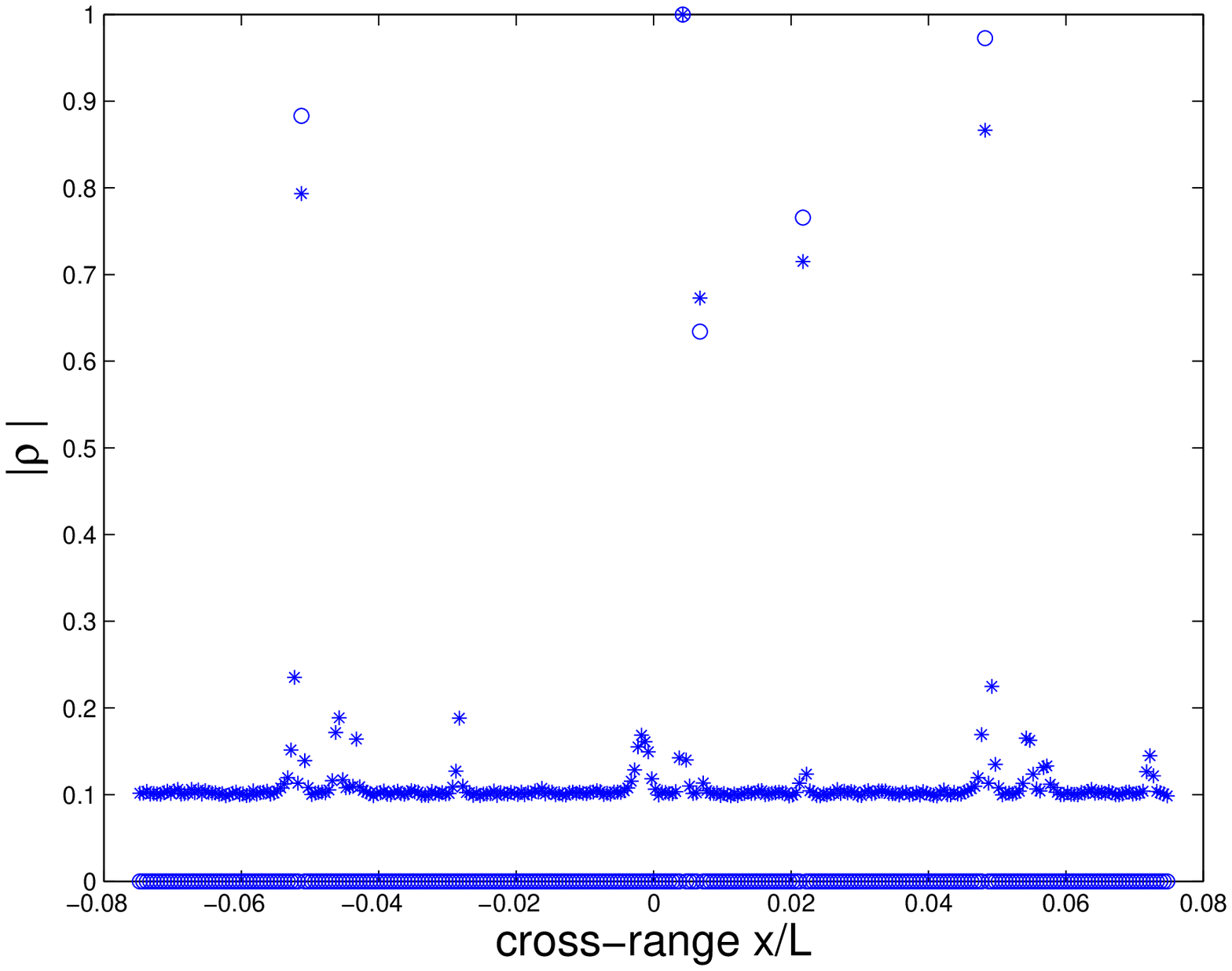} \\
\includegraphics[scale=0.24]{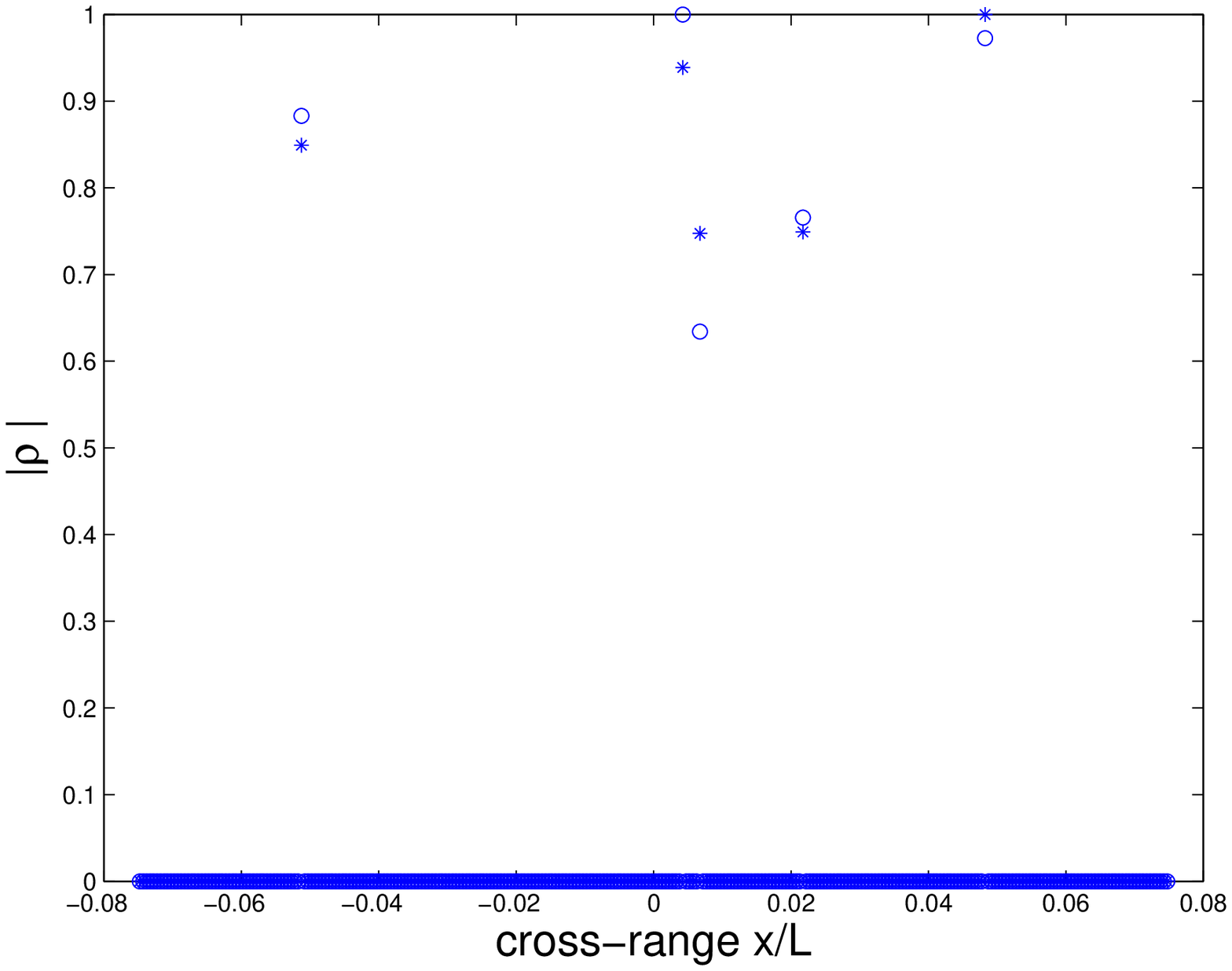} &
\includegraphics[scale=0.24]{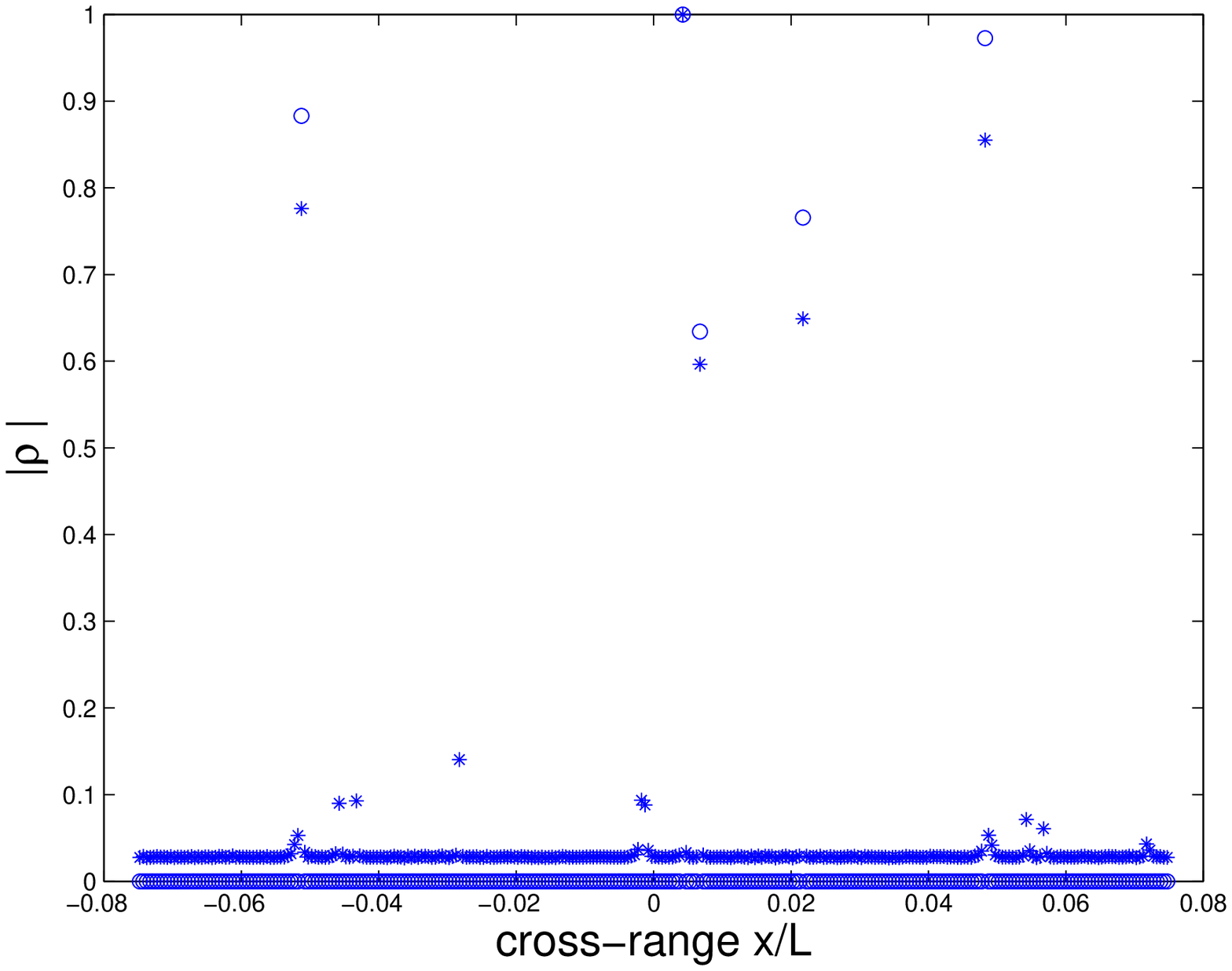} &
\includegraphics[scale=0.24]{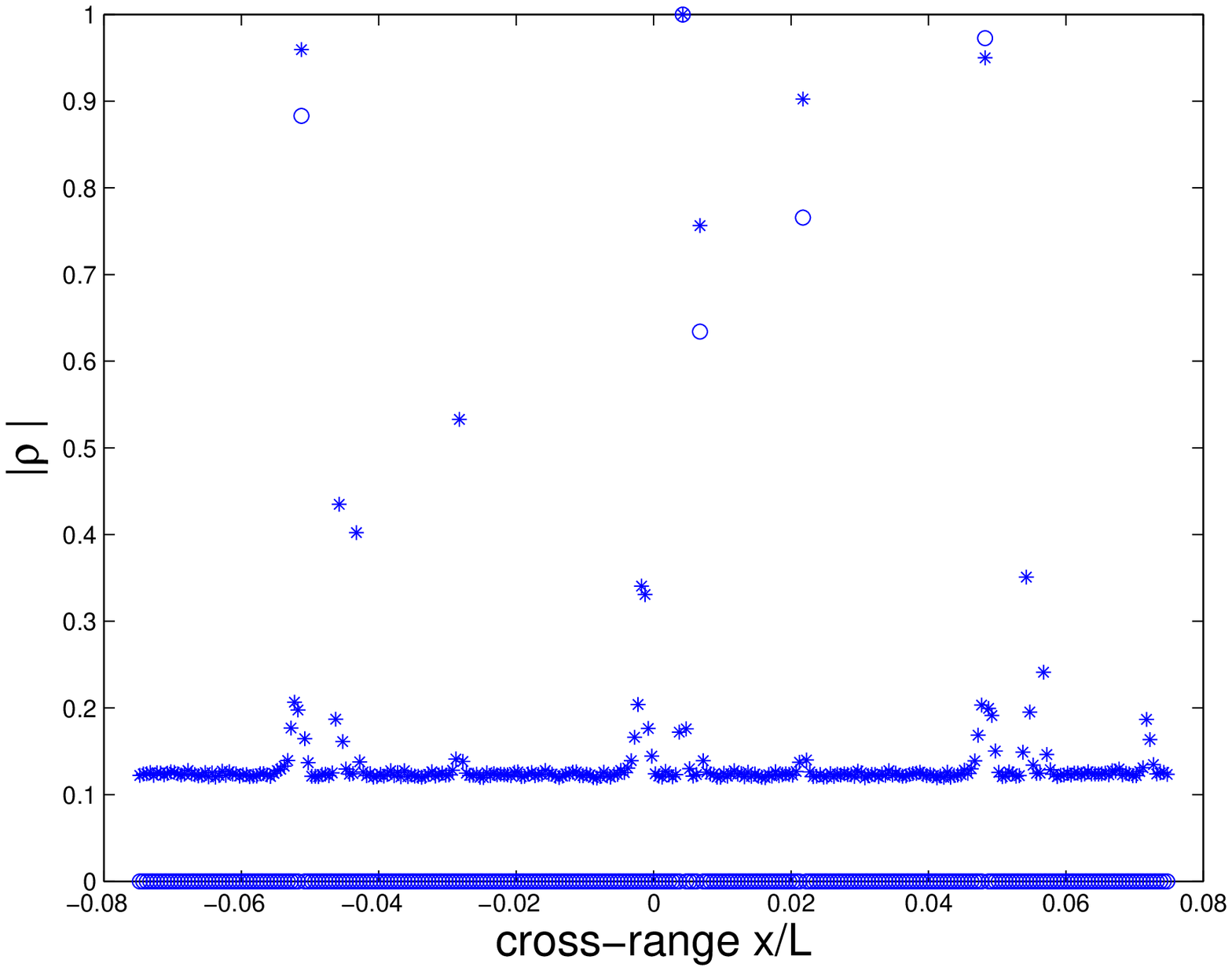} \\
\includegraphics[scale=0.24]{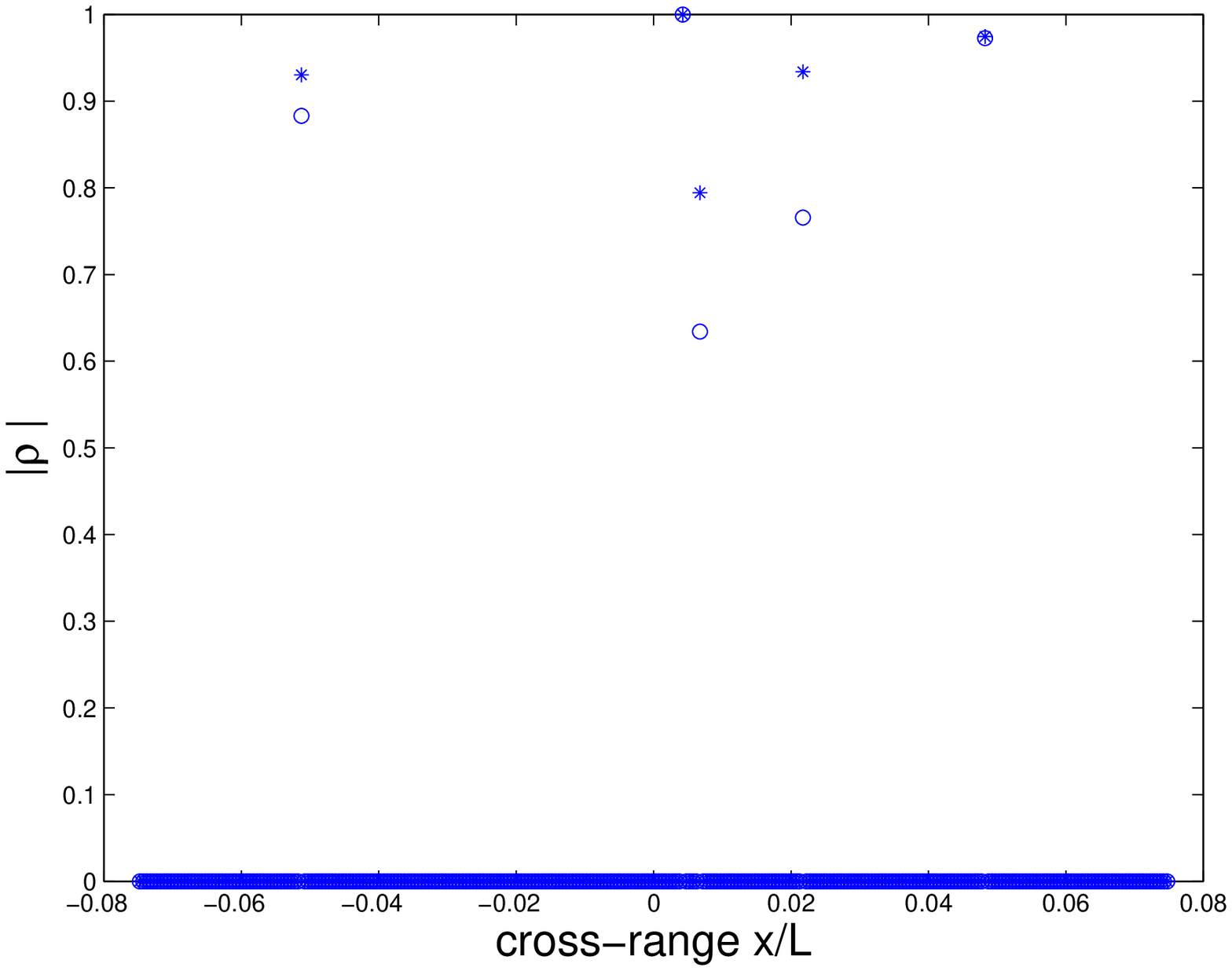} &
\includegraphics[scale=0.24]{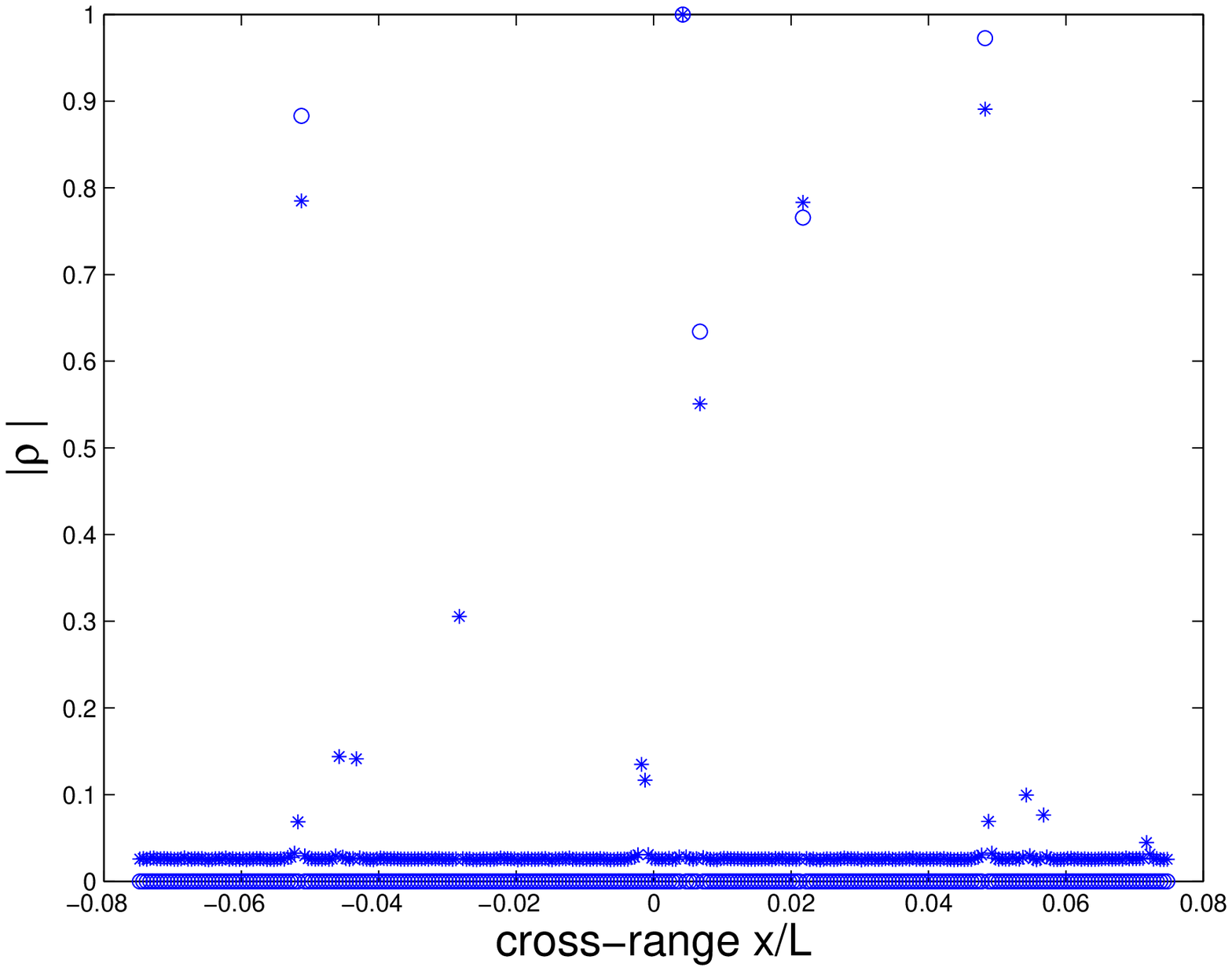} &
\includegraphics[scale=0.24]{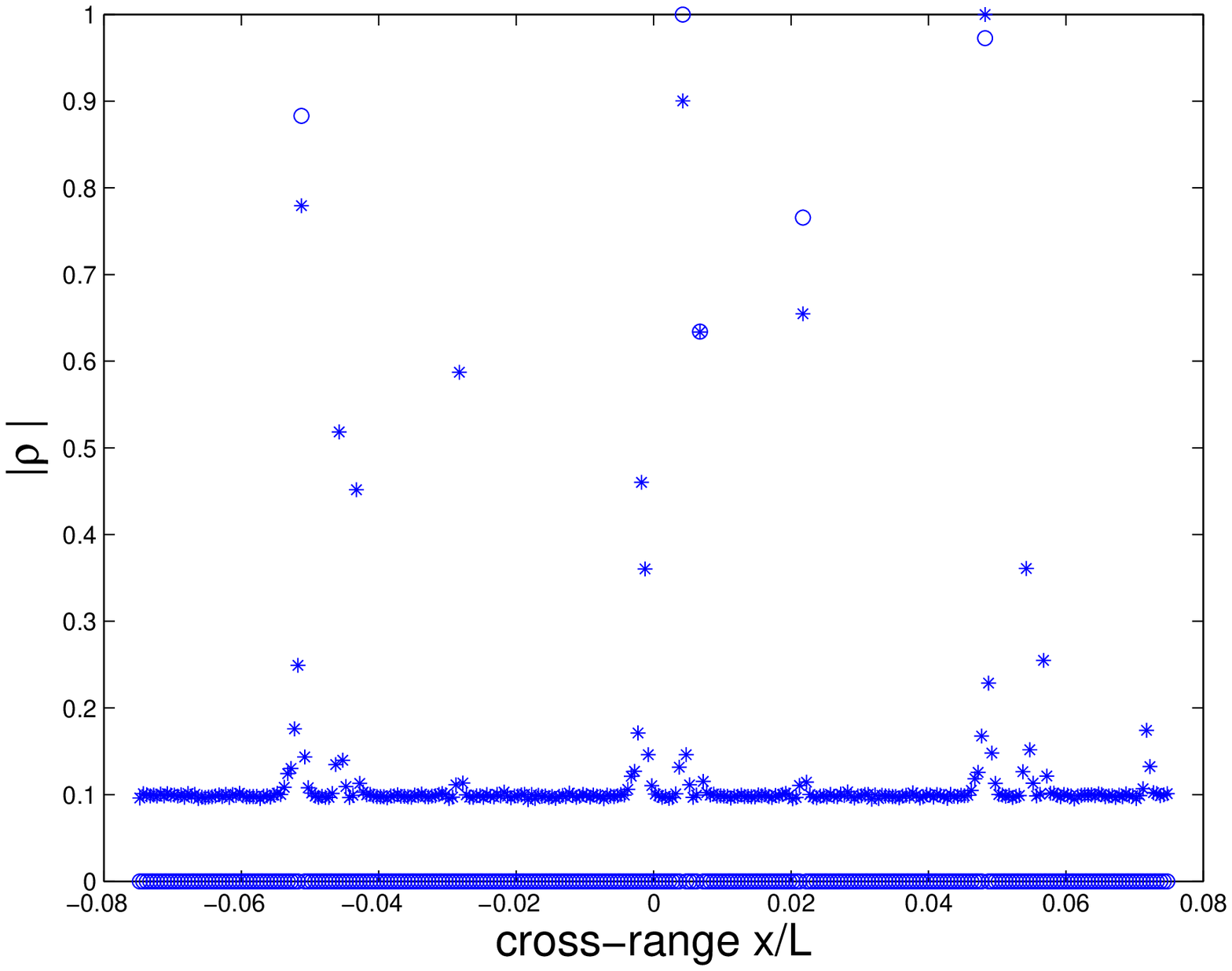} 
\end{tabular}
\caption{Reconstructions obtained with MUSIC from the full recovered matrix $M$. In these numerical simulations
the range $L$ is assumed to be known. 
In the top, middle and bottom rows the ranges are $L=10000\lambda$, $L=20000\lambda$ and $L=50000\lambda$, respectively.
The left, middle and right columns contain $0\% $,  $10\% $, and $30\% $ of additive noise in the data, respectively. 
Length of the array $a=2000 \lambda$, and number of transducers $Na=101$. 
The  cross-range axes are in units of $h_x=\lambda L/a$.
\label{fig:norange}
}
\end{figure}

To study the robustness of the proposed approach in the cross-range direction with respect to additive noise we consider flat images in Figure~\ref{fig:norange}. By a flat image we mean that all the scatterers are placed at the
same range $L$ as in Fig. \ref{fig:schematic} (b), which we assume to be known. The unknown positions of the scatererers are $(x_i,0,L)$,
$i=1,\dots,M$. The exact locations of the scatterers in these images are indicated with circles, and the 
peaks of the MUSIC pseudo-spectrum with stars.
In the top row, the scatterers are placed at range $L=10000\lambda$, in the middle one
at $L=20000\lambda$, and in the bottom one at $L=50000\lambda$.
The data used for the results in the left, middle and right columns contain $0\% $,  $10\% $, and $30\% $ of additive noise, respectively. 
Again, we obtain perfect cross-range positions when noiseless data is used. Furthermore, in the top row ($L=10000\lambda$) all the peaks of 
the MUSIC pseudo-spectrum above the noise level correspond to scatterer locations. However, in the middle row ($L=20000\lambda$)  
a few ghosts appear above the noise level when $30\% $ of noise is added the data. 
In the bottom row ($L=50000\lambda$) the ghosts appear above the noise level  even with  $10\%$ of  noise added to data. We conclude, as expected, 
that the robustness of the MUSIC algorithm is affected by the distance between the array and the IW. This is so because the larger the distance the 
flatter are the wavefronts that illuminate the IW.

\begin{figure}[!htb]
\centering
\begin{tabular}{ccc}
\includegraphics[scale=0.24]{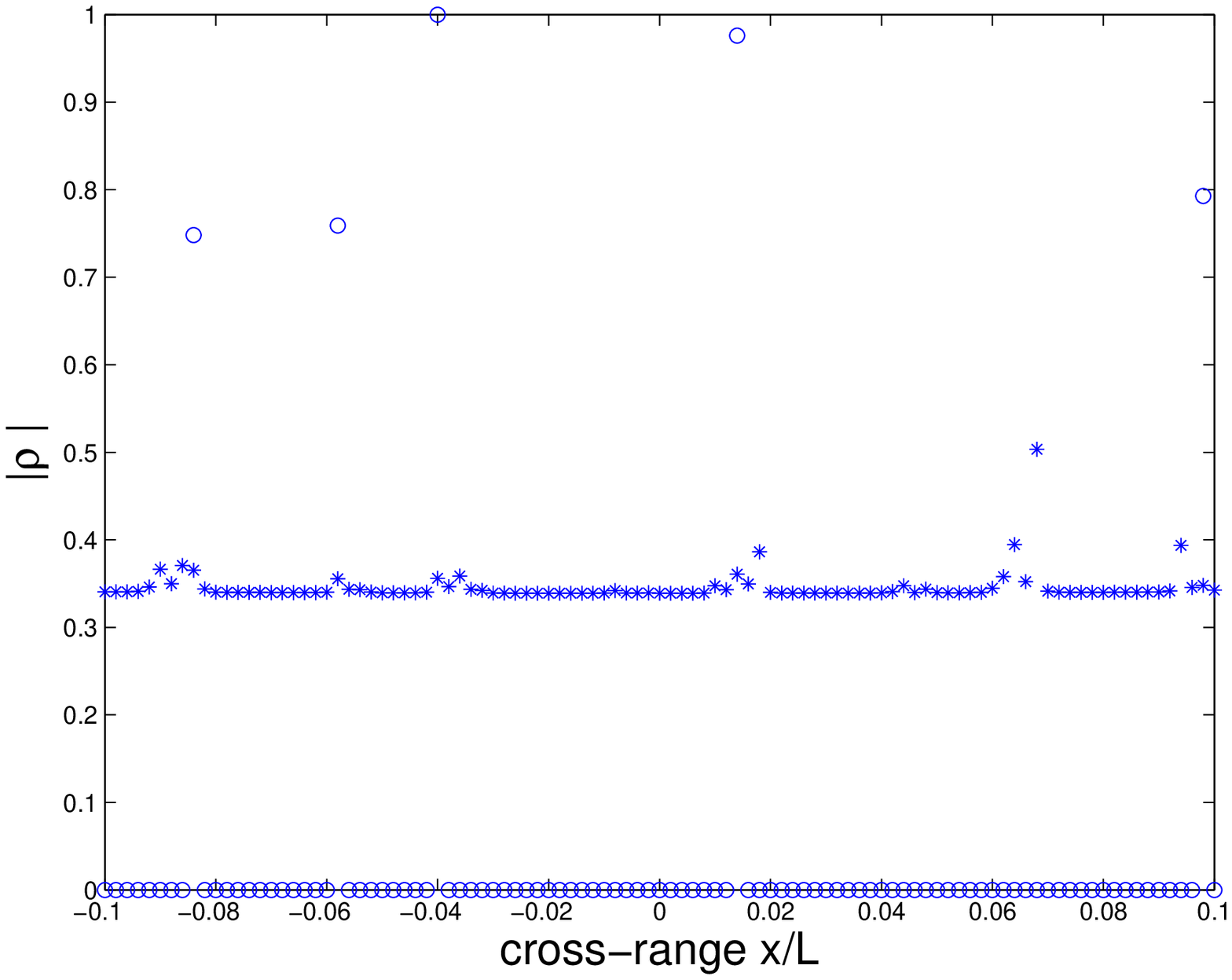} &
\includegraphics[scale=0.24]{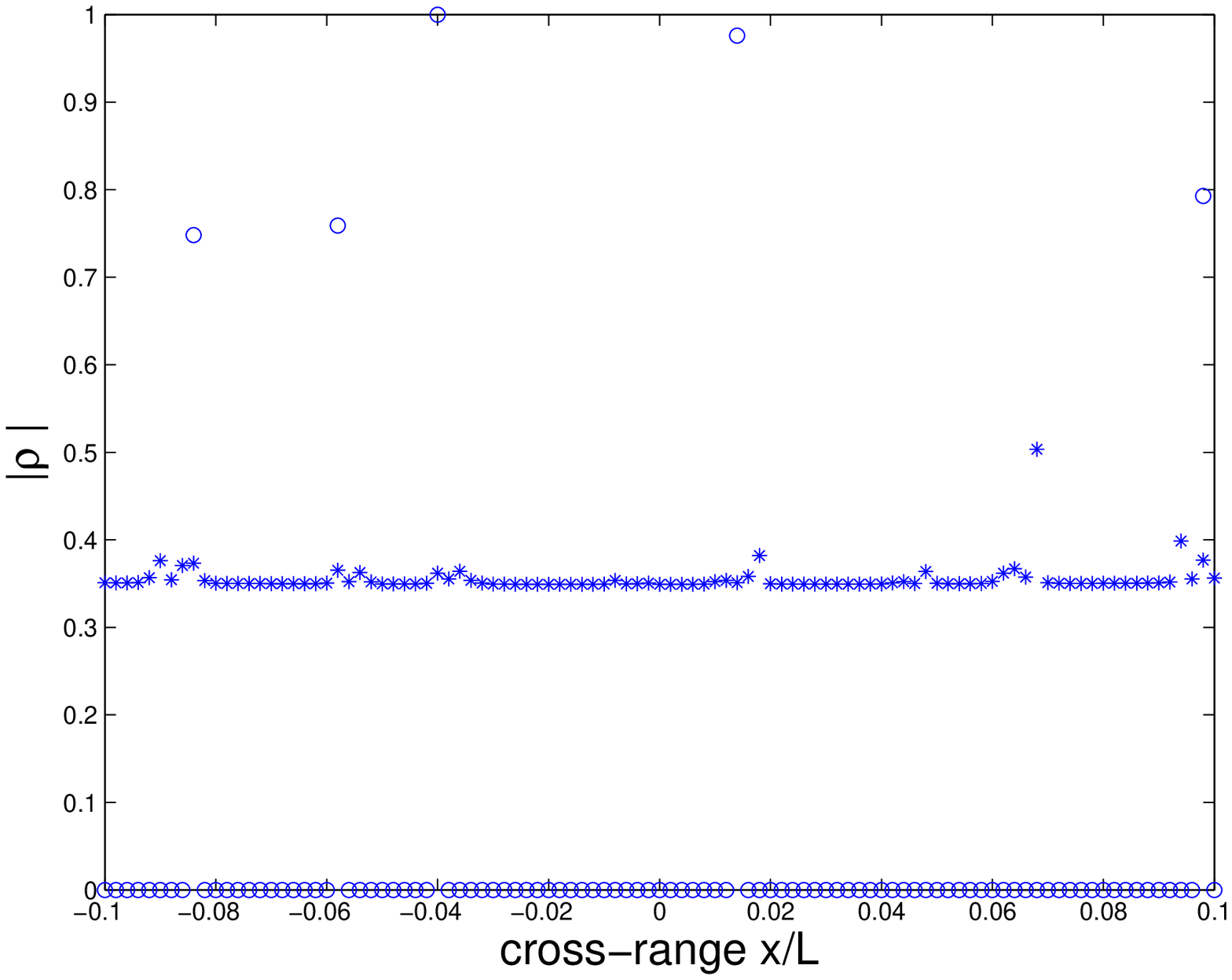} &
\includegraphics[scale=0.24]{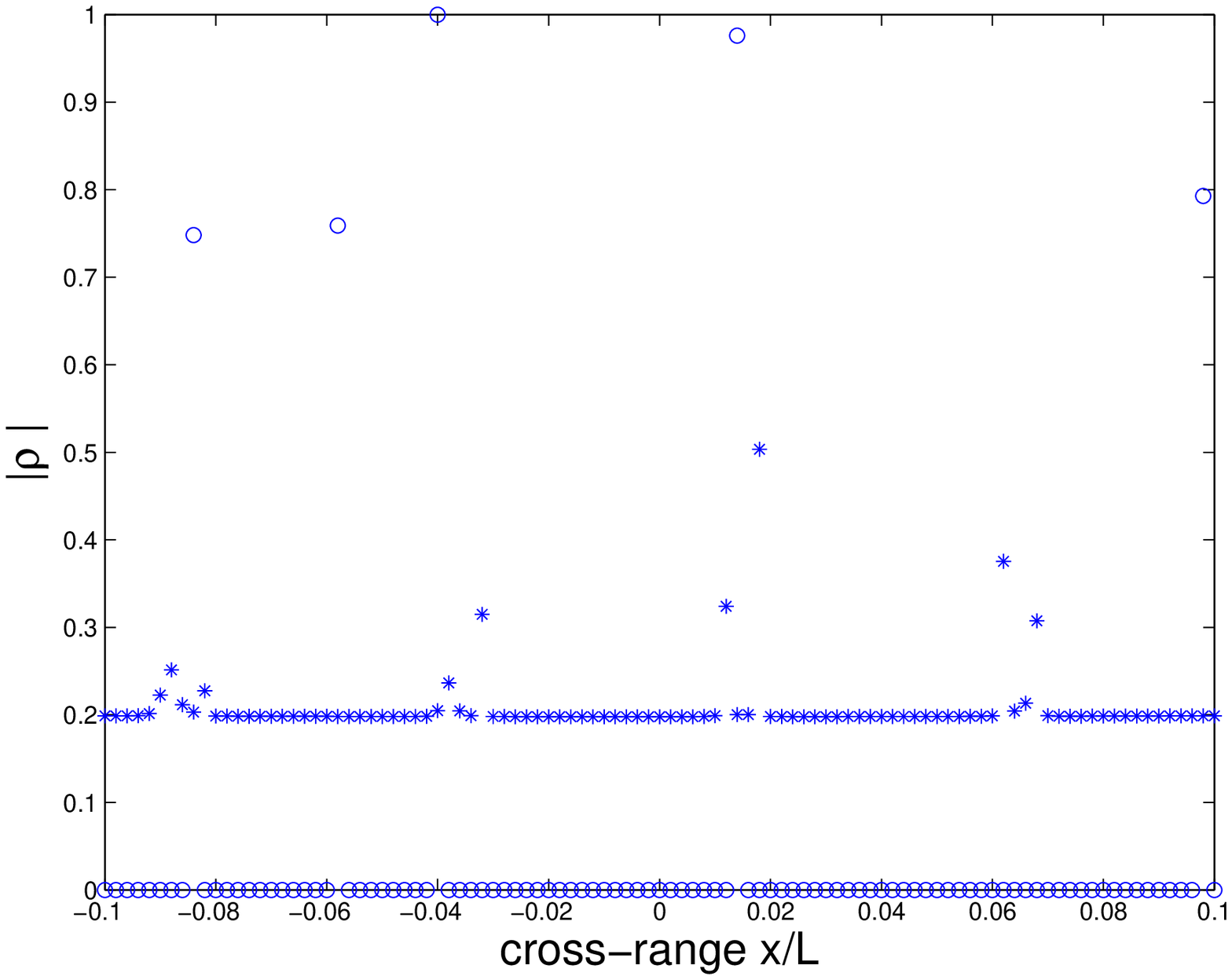} \\
\includegraphics[scale=0.24]{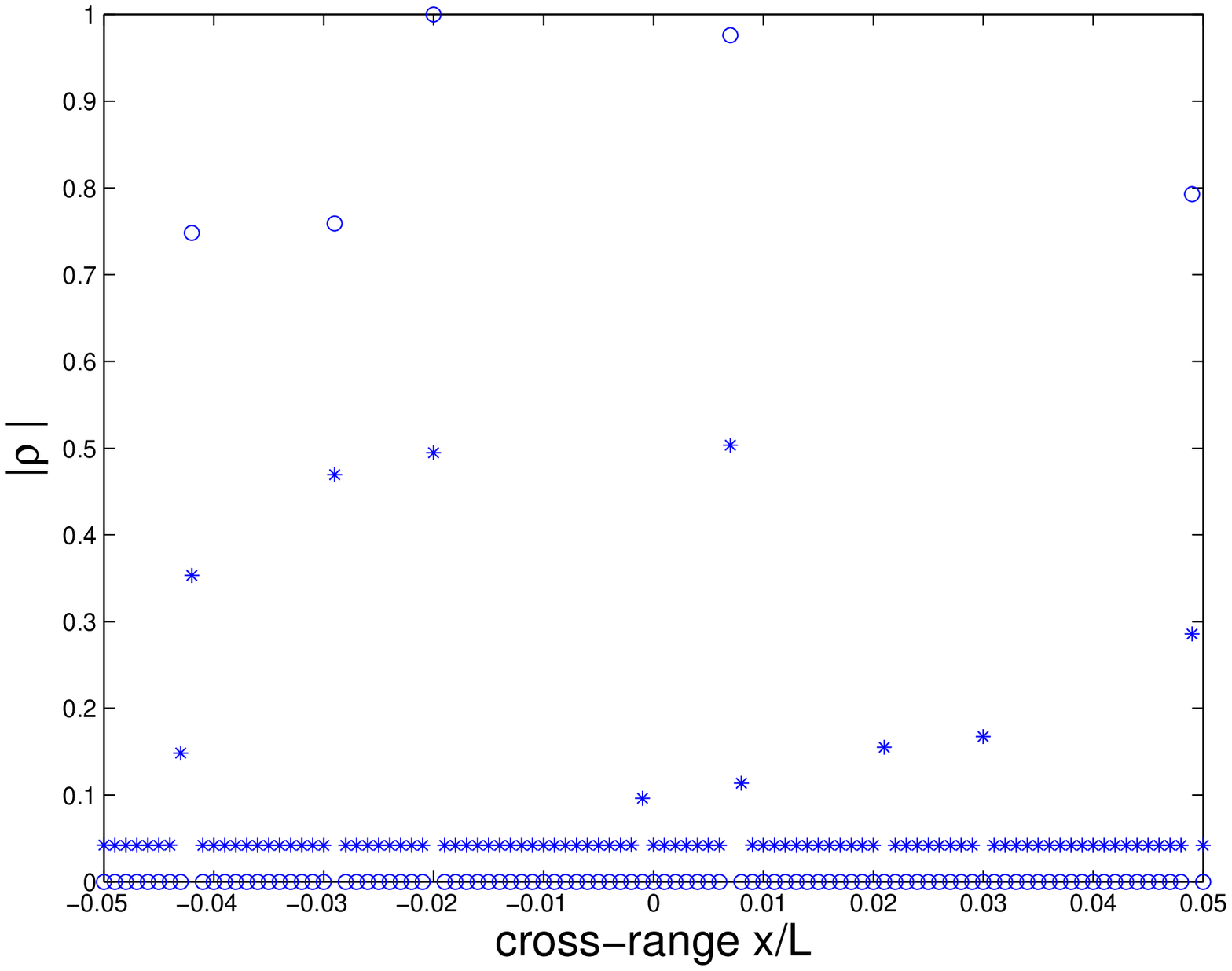} &
\includegraphics[scale=0.24]{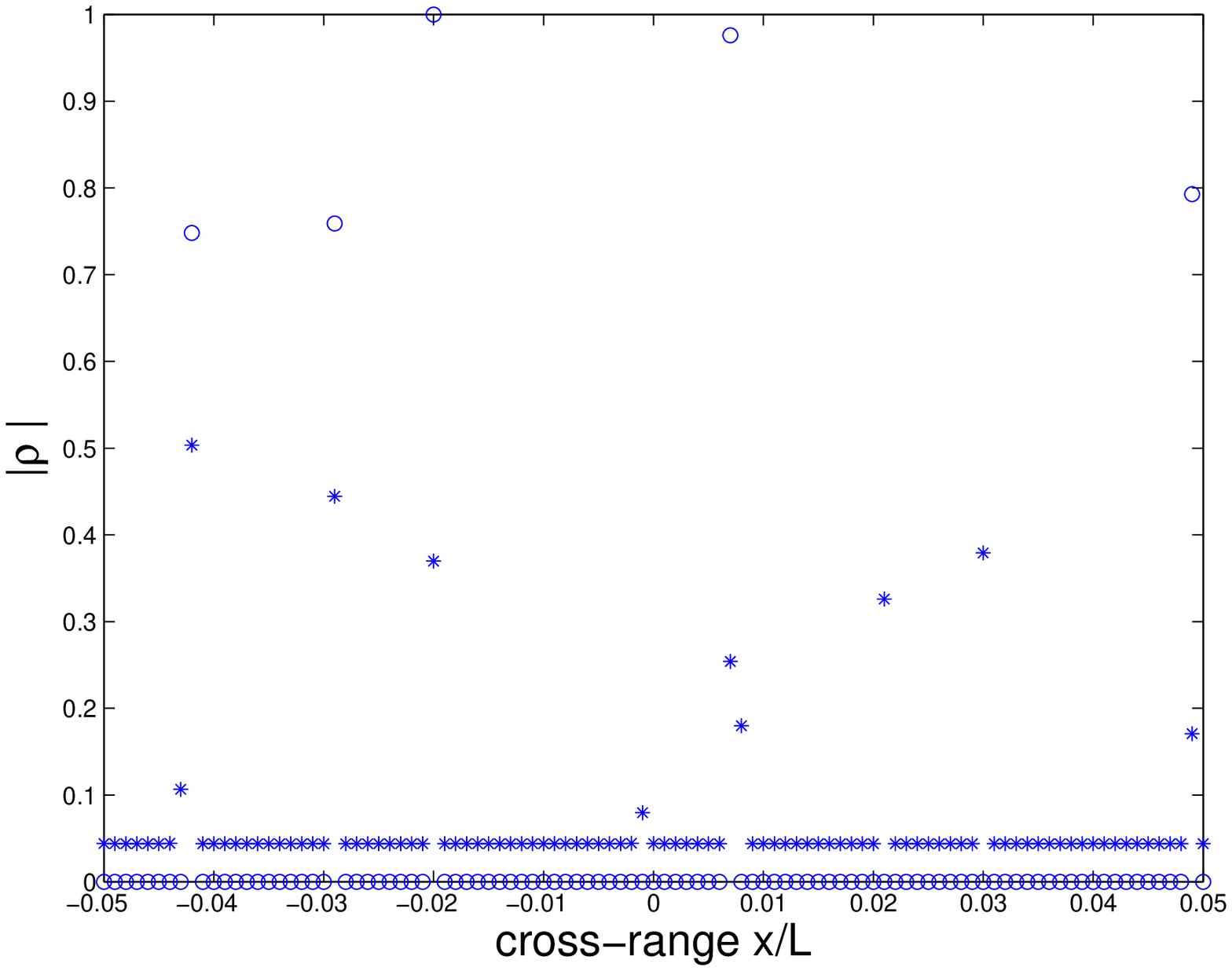} &
\includegraphics[scale=0.24]{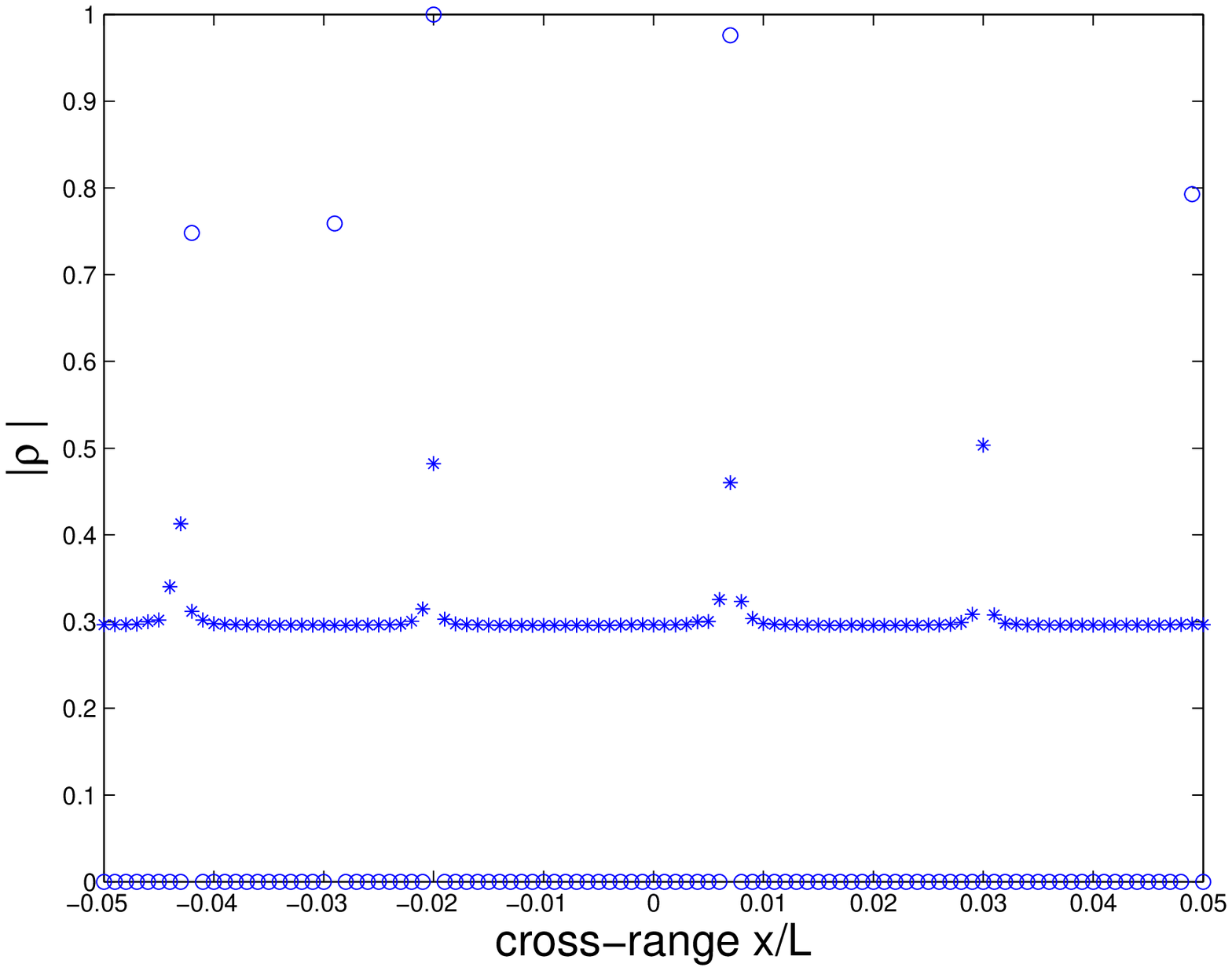} \\
\includegraphics[scale=0.24]{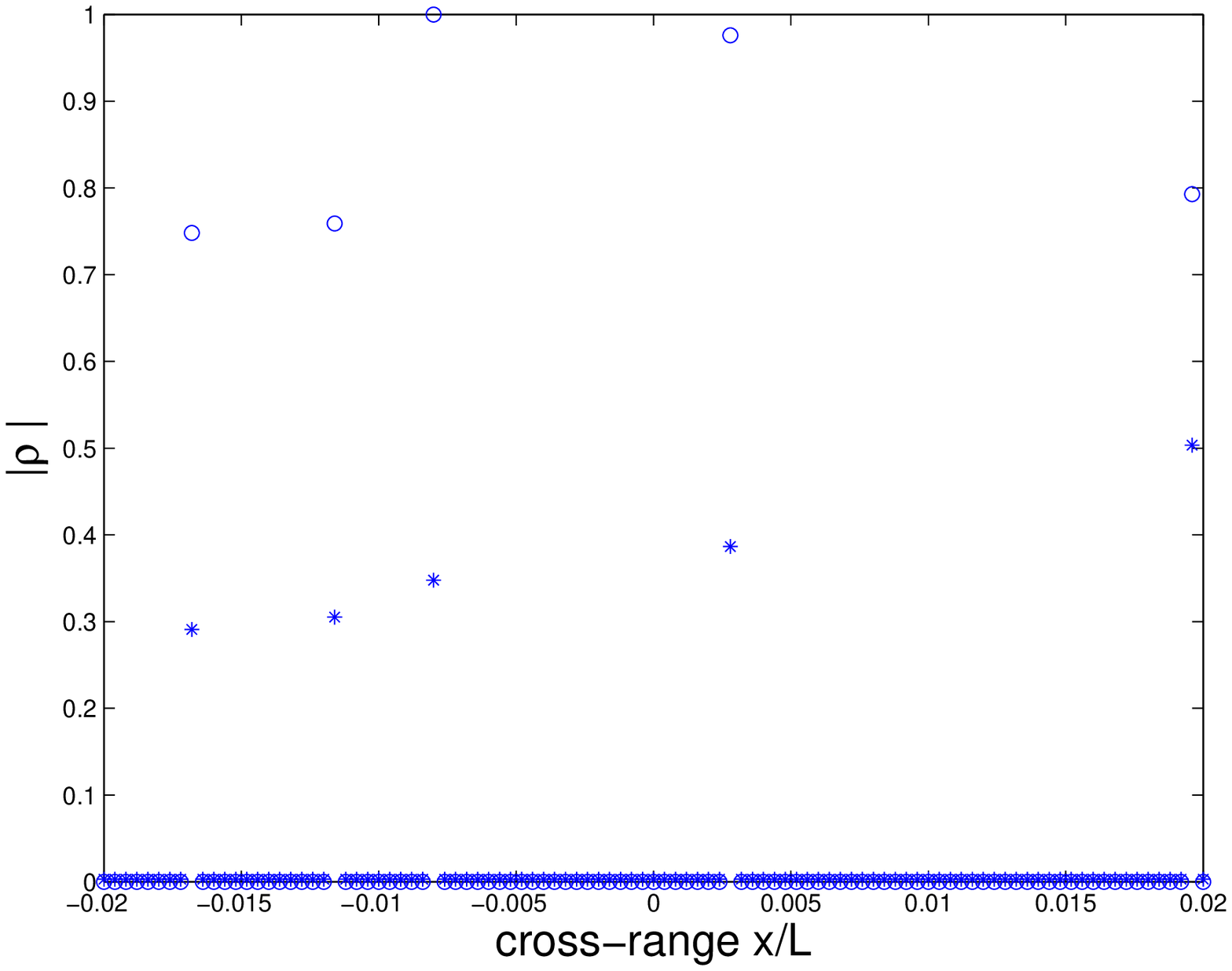} &
\includegraphics[scale=0.24]{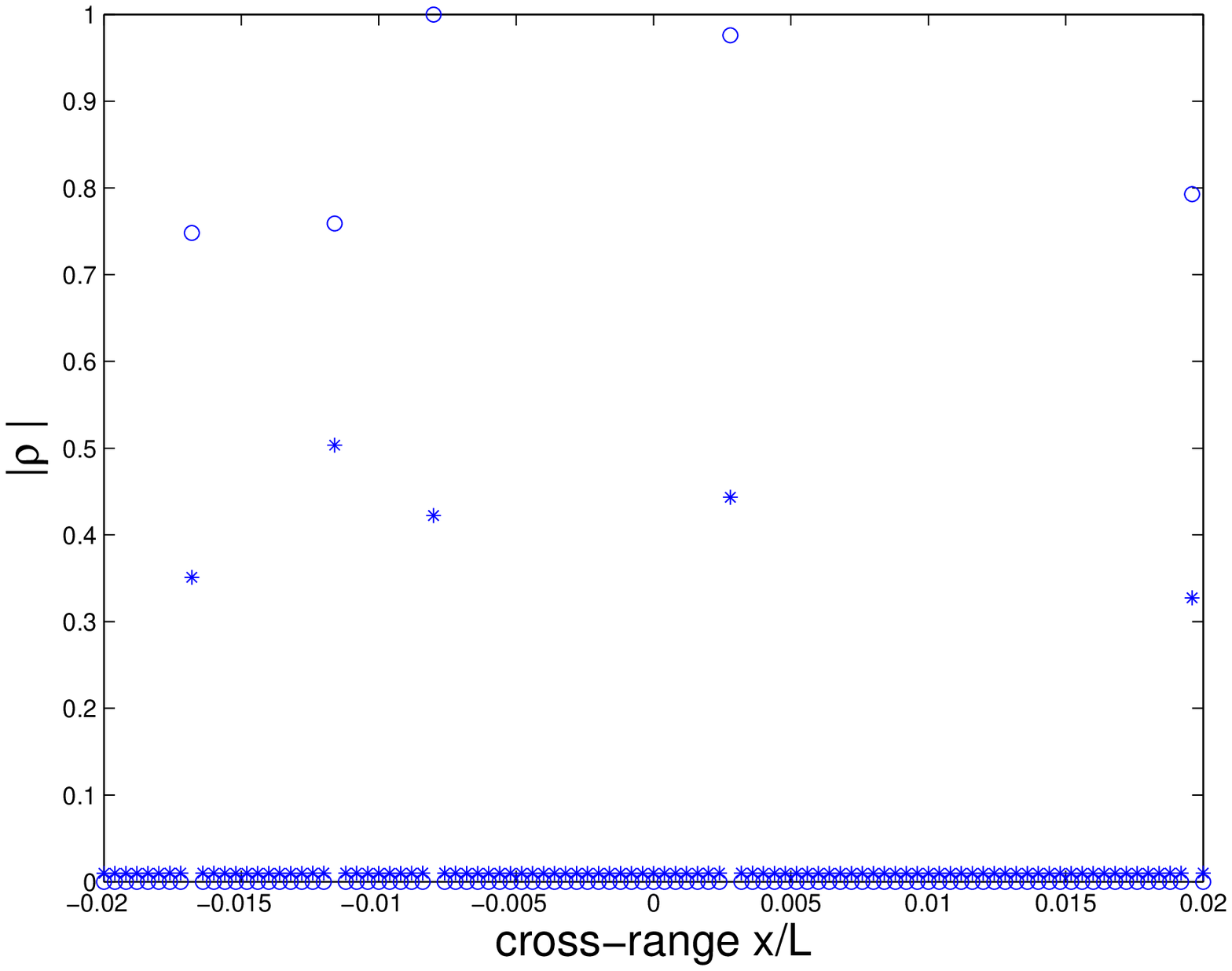} &
\includegraphics[scale=0.24]{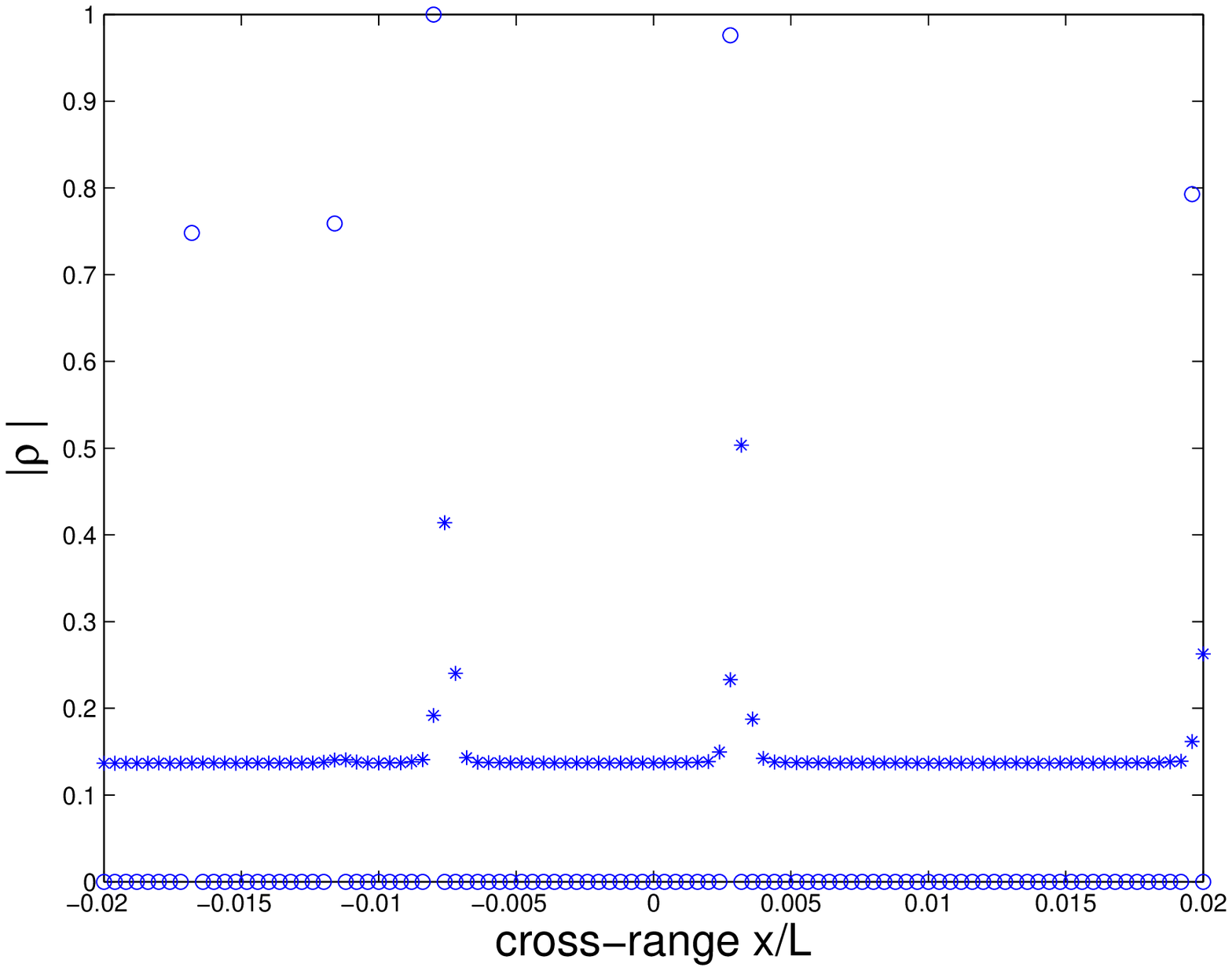}  \\
\includegraphics[scale=0.24]{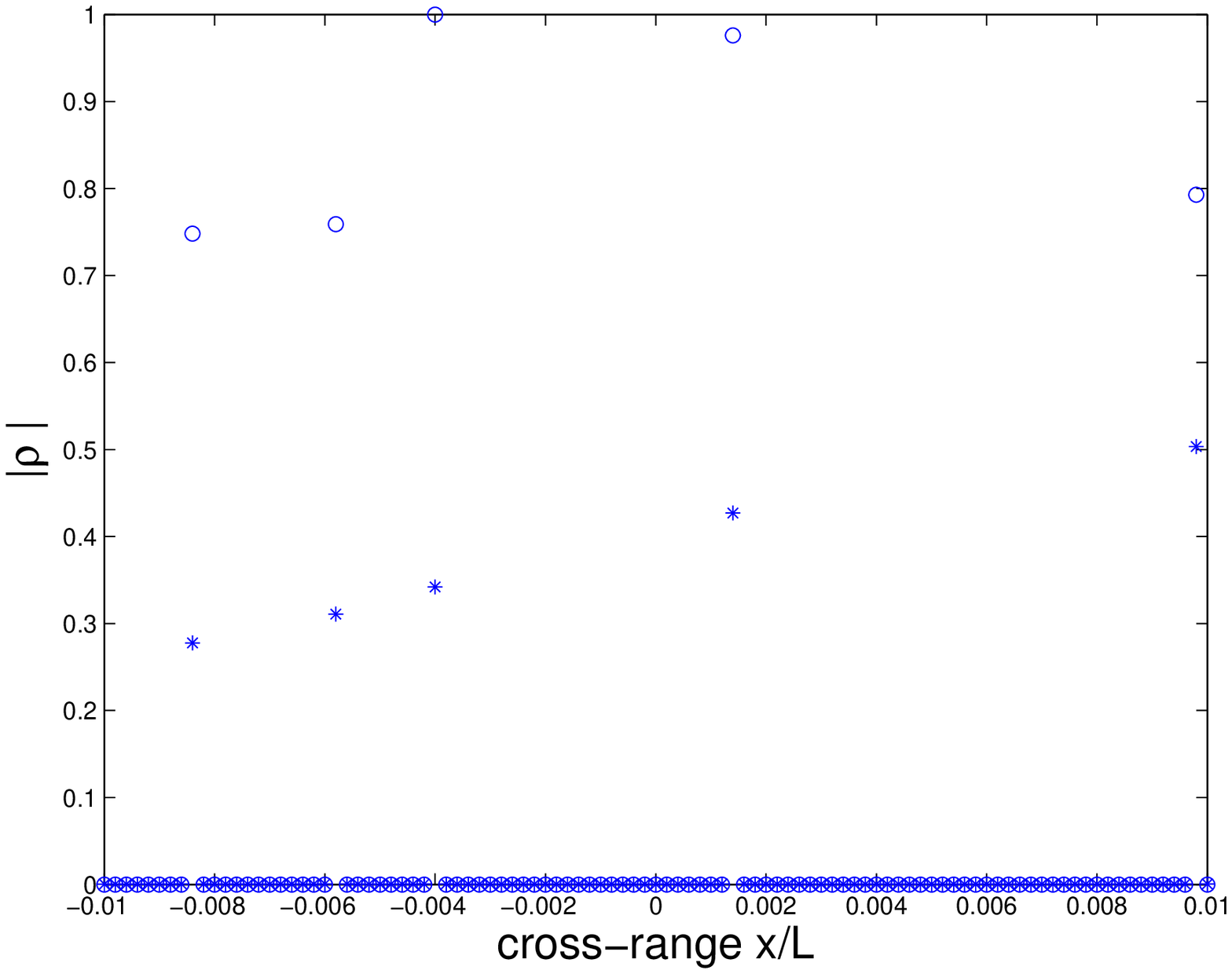} &
\includegraphics[scale=0.24]{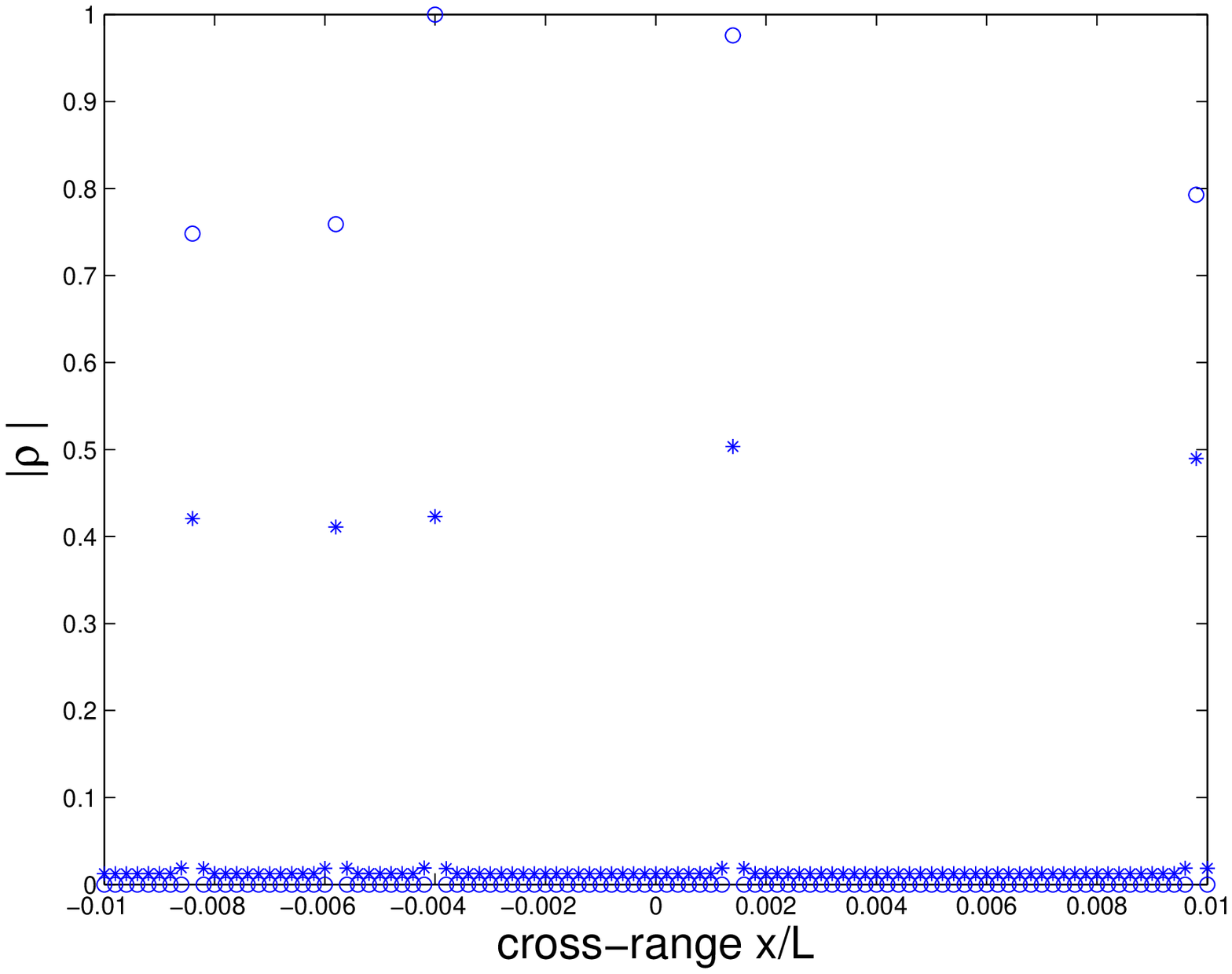} &
\includegraphics[scale=0.24]{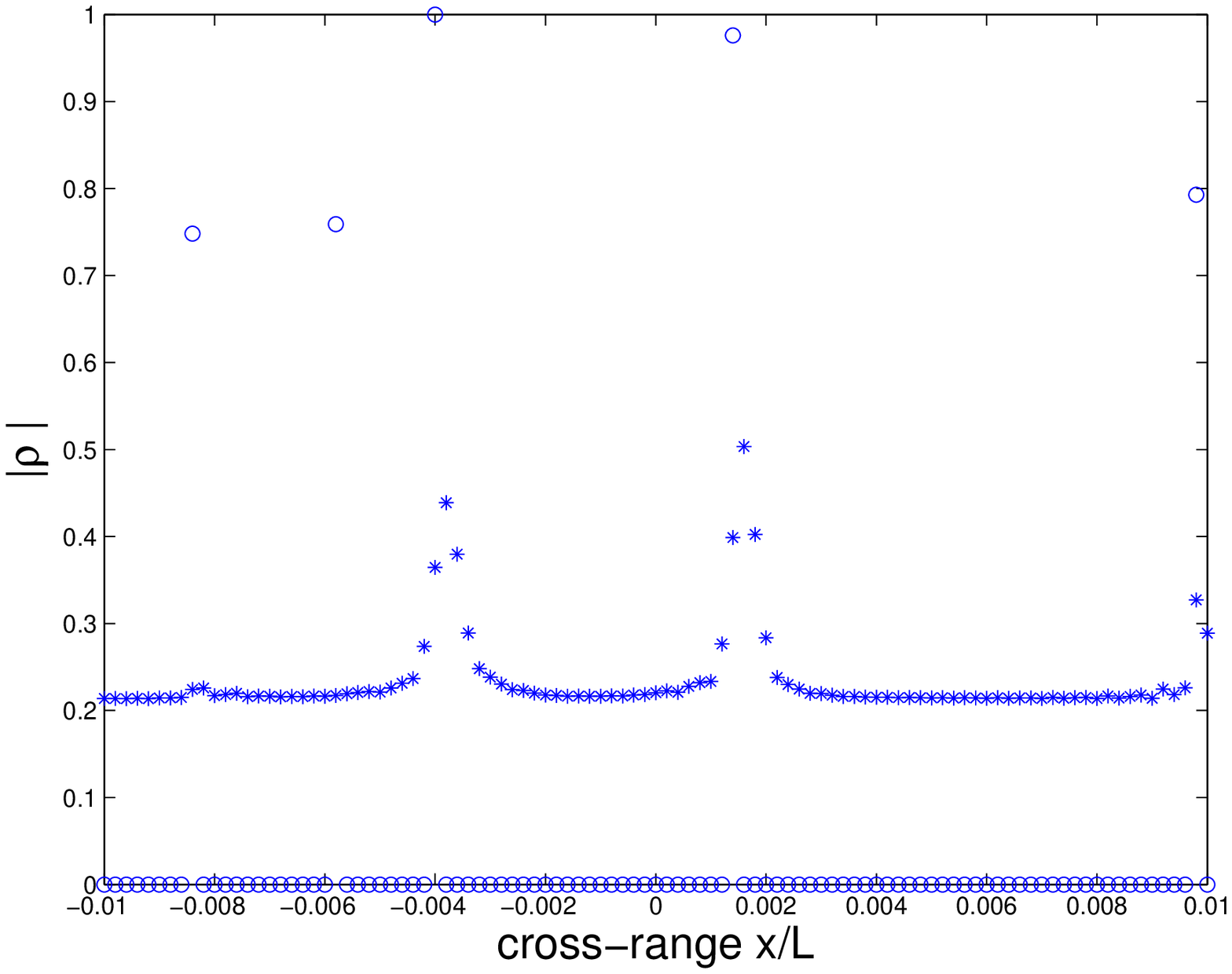}
\end{tabular}
\caption{Reconstructions obtained using MUSIC for the matrix $\vect P$ recovered from six illuminations only. 
In the first, second, third, and forth rows the distances between the array and the IW are $L=10000\lambda$, $L=20000\lambda$on, and $ L=50000\lambda$, 
and $L=100000\lambda$, respectively.
The left, middle and right columns contain $0\% $,  $1\% $, and $10\% $ of additive noise in the data, respectively. 
The length of the array is $a=2000 \lambda$, and number of transducers is $Na=101$. 
The  cross-range axes are in units of $h_x=\lambda L/a$.
\label{fig:hankel}
}
\end{figure}

Our next set of numerical simulations is related to the method of six illuminations.
In Figure~\ref{fig:hankel}, we show the results when wave propagation can be approximated by the paraxial approximation and 
the full response matrix $\vect P$ can be recovered, up to a global phase, using six illuminations as it is explained in subsection~\ref{sec:illum_six}.
In this figure, the scatterers are located on a line at range $L=10000\lambda$ (first row), $L=20000\lambda$ (second row), $L=50000\lambda$ (third row), 
and  $L=100000\lambda$ (fourth row). The exact locations of the scatterers 
are indicated with circles, and the peaks of the MUSIC pseudo-spectrum of the recovered  matrix $\vect P$ with stars. 
The images shown in the left, middle, and right columns are formed using data that contain $0\% $,  $1\% $, and $10\% $ of additive noise, respectively. 
The synthetic data is generated using \eqref{dataI}, with $\vect P$ given in \eqref{responsematrix}, i.e., the data is not generated using the paraxial model 
\eqref{eq:fresnel}.
For noiseless data (left column), Figure~\ref{fig:hankel} shows perfect reconstructions when the distances between the array and the IW are 
large enough (see the fourth and third rows for $L=100000\lambda$ and $L=50000\lambda$, respectively).
However, the reconstructions become poor as we move the IW closer to the array (see the left column of the
second and first rows for $L=20000\lambda$ and $L=10000\lambda$, respectively). This is due to modeling
errors as the paraxial approximation deteriorates as we move the IW closer to the array.
Because the distances between the array and the IW are very large, we also observe that the method is less robust with respect to additive noise. 
With $10\%$ of noise the method fails to locate all the scatterers for all $L$ (see the right column). Only for 
low level of noise and large enough $L$ the method locates all the scatterers (see the middle column).

\section{Conclusions}
\label{sec:conclusions}

In this paper we consider narrow band, active array imaging of weak localized scatterers
when only the intensities are recorded and measured at an array with $N$ transducers.  We assume that the medium is homogeneous so wave 
propagation is coherent.
We show that if one controls the illuminations, the imaging problem with intensity-only can be easily reduced to one in which the phases are available. 
Furthermore, the amount of extra work in data acquisition is small, as  only $3N-2$ illuminations are needed in general.  
We also show that Fresnel and Fraunhofer diffraction images can be obtained by using only $6$ illuminations for 1-D images and $12$ illuminations
for 2-D images, if these images are flat and their ranges are known. 

Numerical simulations show the performance of the proposed approaches with and without additive noise in the data.
The results presented in the paper indicate different noise sensitivities of these approaches. In particular, the method of six illuminations, that is valid when the distance between the scatterers and the 
array is much larger than the wavelength of the signals and much larger than the linear dimensions of the array and the IW, is very sensitive to 
noise. This limits the use of the six illuminations approach in 
practice to situations where the signal to noise ratio is high. Detailed analysis of robustness to noise of the proposed 
phaseless imaging methods will be presented in our consequent work. 

\section*{Acknowledgment} 

Miguel Moscoso's work was partially supported by  AFOSR FA9550-14-1-0275 and the Spanish MICINN grant
FIS2013-41802-R. George Papanicolaou's work was partially supported by
AFOSR grant FA9550-14-1-0275. Alexei Novikov's work was partially supported by  AFOSR FA9550-14-1-0275 and NSF DMS-1515187.

\appendix
\section{Derivation of the paraxial approximation}\label{seq:app}
Let the array be located on the plane $\vec{\vect x}=(\vect x,0)$ and the scatterers on the plane 
$\vec{\vect y}=(\vect y,L)$. Then, we have
\begin{equation}
 \abs{\vec{\vect x} - \vec{\vect y}}=L\left(1 + \frac{1}{2}\left(\frac{ \abs{\vect x - \vect y}}{L}\right)^2 + 
O\left(\frac{ \abs{\vect x - \vect y}}{L}\right)^4 \right)\, .
\label{eq:taylor}
\end{equation}
If the distance $L$ between the center of the array and the center of the IW  is  large enough
and the linear dimensions of the array and the IW are small, so only the two first terms in \eqref{eq:taylor} 
have to be retained, then \eqref{response} becomes 
\begin{equation}
P(\vect x_r,\vect x_s) \approx  
\Pp(\vect x_r,\vect x_s) = \frac{e^{i2 \kappa L}}{(4\pi L)^2} 
\sum_{j=1}^M \alpha_j e^{i \kappa \left[\left| \vect x_s - \vect \xi_j \right|^2+ \left|
 \vect x_r - \vect \xi_j \right|^2 \right]/2L} \,,
\label{eq:parax0}
\end{equation}
where ${\vect \xi}_j$ are the two-dimensional cross-range vectors that denote the scatterer's  positions, and
$\vect x_s$ and $\vect x_r$ are  the two-dimensional cross-range vectors that denote the positions of the sources and the receivers, 
respectively. If we further expand the terms in the exponential, we obtain
\begin{equation}\everymath{\displaystyle}
\Pp(\vect x_r,\vect x_s)  = C_r C_s  \sum_{j=1}^M \bar{\alpha}_j e^{- i \kappa \frac{ \langle \vect x_s + \vect x_r,  
\vect \xi_j \rangle}{L}}  = C_r C_s \sum_{j=1}^K \tilde{\rho}_j e^{- i \kappa \frac{\langle \vect x_s + \vect x_r,  
\vect \xi_j \rangle}{L}}  \,,
\label{eq:fresnel_o}
\end{equation}
where 
\begin{equation}\everymath{\displaystyle}
\tilde{\rho}_j = \rho_j e^{i \kappa \left| \vect \xi_j \right|^2/2L}\, 
\end{equation}
are distorted reflectivities that only change the phases of the reflectivities of the scatterers, and
\begin{equation}\label{eq:factor_o}
 C_t =C(\vect x_t) = \frac{e^{i  \kappa L}e^{i \kappa \left| \vect x_t \right|^2 /2L}}{4\pi L},~t=r,s,
\end{equation}
are geometric factors that only depend on the imaging setup.

Equation \eqref{eq:fresnel_o} gives accurate results if the approximation in the phase for any two points $\vec{\vect x}$ and $\vec{\vect y}$ in the array and the IW is such that contribution of the third term in \eqref{eq:taylor} is
small, i.e., if
\begin{equation}
\frac{\kappa\abs{\vect x - \vect y}^4}{8 L^3}= \frac{\pi\abs{\vect x - \vect y}^2}{4\lambda L}\left( \frac{\abs{\vect x - \vect y}}{L}\right)^2 \ll\pi\, .
\label{eq:thirdterm}
\end{equation}
Let $a$ and $b$ be two characteristic lengths of the array and the IW, respectively. 
Typically $a\leq b$, so
$\abs{\vect x - \vect y} \leq a$ for all $\vect x$ and $\vect y$. 
Therefore, \eqref{eq:fresnel_o} models of wave propagation accurately if 
\begin{equation}
\frac{1}{4} F \left(\frac{a}{L}\right)^2 \ll 1\, ,
\label{eq:cond1Fresnel}
\end{equation}
where  
\begin{equation}
F = \frac{a^2}{\lambda L}\, 
\label{eq:fresnelnumber0}
\end{equation}
is the Fresnel number.
In \eqref{eq:parax0} and \eqref{eq:fresnel_o} we have kept the quadratic term of the expansion \eqref{eq:taylor}. 
In order for this term to be significant, we must have that $\frac{\kappa}{2}\frac{a^2}{L}\gtrsim\pi$, which means that
\begin{equation}
F\gtrsim 1\, 
\label{eq:cond2Fresnel}
\end{equation}
must hold. Conditions \eqref{eq:cond1Fresnel} and \eqref{eq:cond2Fresnel} characterize the Fresnel or near field regime. In this regime, the array does not have to be too small so the wave fronts appear to be planar when viewed from the array, as is the case in the Fraunhofer or far field regime. 
In the Fraunhofer regime we have the condition $F\ll 1$ instead, which means that the quadratic phase terms in \eqref{eq:taylor} are negligible.

\section{MUSIC}\label{Be}
MUSIC is a subspace projection algorithm that uses the 
SVD of the full data array response matrix $\vect\wP(\omega)$ to form the images. 
It is a direct algorithm widely used to image the locations of $M<N$ point-like scatterers in a region of interest. 
Once the locations are known, their reflectivities can be found from the recorded intensities using convex optimization as shown below.

Let us assume in this Section that $\vect\wP(\omega)$ is fully recorded and known. 
We write the SVD of the data matrix $\vect\wP(\omega)$ in the form
\begin{equation}
\label{svd1}
\vect\wP(\omega)=\vect\wU(\omega)\vect\Sigma(\omega)\vect\wV^\ast(\omega)
=\sum_{j=1}^{\tilde M}\sigma_j(\omega)\wU_j(\omega)\wV_j^\ast(\omega)\, ,
\end{equation}
where $\sigma_1(\omega)\ge\cdots\ge\sigma_{\tilde M}(\omega)>0$ are the nonzero singular values, 
and $\wU_j(\omega)$, $\wV_j(\omega)$ are the corresponding left and right singular vectors, respectively. They fulfill the following equations:
\begin{equation}
\label{svd2}
\vect\wP^\ast(\omega) \wU_j(\omega) = \sigma_j(\omega) \wV_j(\omega)\, , \quad 
\vect\wP(\omega) \wV_j(\omega) = \sigma_j(\omega) \wU_j(\omega)\,,\,\,j=1,\ldots,N .
\end{equation}
Since $\vect\wP(\omega)$ is symmetric, $\wU_j(\omega)= e^{{\bf i} \theta_j} {\overline \wV_j(\omega)}$  for some  unknown global phase $\theta_j$, $j=1,\ldots,N$.

The search of the locations  of the $M$ scatterers is the combinatorial part of the imaging problem and, hence, by far the most difficult task. 
Note that $\vect\wP(\omega)$ is a linear transformation from the {\em illumination space} $\mC^N$ to the {\em data space} $\mC^N$. According to \eqref{svd1}, the
 illumination space can be decomposed into the direct sum of
a  signal space, spanned by the principal singular vectors $\wV_j(\omega)$, $j=1\dots,M$, having non-zero singular values,
and a noise space spanned by the singular vectors having zero singular values.
Since the singular vectors
$\widehat{V}_j(\omega)$, $j=M+1,\ldots,N$, span the noise space, the probing vectors $\vect\wg_0(\vect y_j,\omega)$ will be orthogonal to the noise space
only when $\vect y_j$ corresponds to a scatterer's location $\vect y_{n_j}$. Hence, it follows that the scatterers' locations must correspond to the peaks
of the functional
\begin{equation}
\label{MUSIC_0}
\mathcal{I}(\vect y_s)=\frac{1}{\sum_{j=M+1}^{N} |\vect\wg_0^T(\vect y_s,\omega)\widehat{V}_j(\omega) |^2 },\,\,s=1,\ldots,K.
\end{equation}
We can interpret \eqref{MUSIC_0} in terms of the images created by the singular vectors having zero singular value, as 
$\vect\wg_0^T(\vect y_s,\omega)\widehat{V}_j(\omega)$ is the incident field at the search point $\vect y_s$ due to a illumination vector 
$\widehat{V}_j(\omega)$ on the array. 
According to this interpretation, the singular vectors having zero singular value do not illuminate the scatterers locations and, hence, \eqref{MUSIC_0}
has a peak when $\vect y_s=\vect y_{n_j}$.

Since in our application the number of scatterers is small, the signal space is much smaller than the noise space and, therefore, it is more efficient to compute
the equivalent functional
\begin{equation}
\label{MUSIC}
\mathcal{I}_{MUSIC}(\vect y_s)=\frac{\min_{1\le j\le K}\|\mathcal{P}\vect\wg_0(\vect y_j,\omega)\|_{\ell_2}}{\|\mathcal{P}\vect\wg_0(\vect y_s,\omega)\|_{\ell_2}},\,\,s=1,\ldots,K,
\end{equation}
with the
projection onto the noise space defined as
\begin{equation}
\label{proyection}
\mathcal{P}\vect\wg_0(\vect y,\omega)=\vect\wg_0(\vect y,\omega) -\sum_{j=1}^M (\vect\wg_0^T(\vect y,\omega)\widehat{V}_j(\omega))
\widehat{V}_j(\omega).
\end{equation}
The numerator in \eqref{MUSIC}  is just a normalization. We note that \eqref{MUSIC}  is robust to noise, even for single frequency and for non-homogeneous, random media, and it is quite accurate
for large arrays~\cite{BTPB02}. Generalizations of MUSIC for multiple scattering and extended scatterers have also 
been developed (see, for example, \cite{Gruber04} and \cite{hou06}).


\end{document}